\documentclass[12pt,a4paper]{article}
\usepackage{latexsym,amssymb,amsmath,mathrsfs}

\usepackage{sectsty}
\usepackage{color}
\usepackage{bm}
\usepackage[T1]{fontenc}
\usepackage[anchorcolor=blue,%
 bookmarks=true,%
 bookmarksnumbered=true,%
]{hyperref}
\usepackage{pst-all}
\usepackage{graphicx}

\allsectionsfont{\normalsize\sc\center}

\newtheorem{thm}{Theorem}[section]
\newtheorem{prop}[thm]{Proposition}
\newtheorem{cor}[thm]{Corollary}
\newtheorem{lem}[thm]{Lemma}
\newtheorem{defn}[thm]{Definition}
\newtheorem{remark}[thm]{Remark}

\newenvironment{pf}{\par\begin{trivlist}%
\item[]{\bf Proof.}\ }{\hfill $\square$ \end{trivlist}\par}
\newenvironment{apf}[1]{\par\begin{trivlist}%
\item[]{\bf Proof of #1.}\ }{\hfill $\square$ \end{trivlist}\par}

\makeatletter \@addtoreset{equation}{section} \makeatother

\newcommand{\R}{\mathbb{R}}

\DeclareMathOperator{\Ric}{Ric}

\renewcommand{\d}{\mathrm{d}}
\newcommand{\m}{\mathfrak{m} }

\newcommand{\s}{\mathfrak{s}}
\newcommand{\vol}{\rm{vol}}

\newcommand{\eps}{{\varepsilon}}

\newcommand{\wh}{\widehat}

\title{\large\bf Lower \boldmath$N$-weighted Ricci curvature bound with \boldmath$\varepsilon$-range and 
displacement convexity of entropies}
\author{Kazuhiro Kuwae\thanks{Department of Applied Mathematics, Fukuoka University,
Fukuoka 814-0180, Japan ({\sf kuwae@}  {\sf fukuoka-u.ac.jp}) . Supported in part by JSPS Grant-in-Aid for Scientific Research (KAKENHI) 17H02846 and by fund (No.:185001) from the Central Research Institute of Fukuoka University.}
\ \ and\ \ 
Yohei Sakurai\thanks{Department of Mathematics, Saitama University,
255 Shimo-Okubo, Sakura-ku, Saitama-City, Saitama, 
338-8570, Japan ({\sf ysakurai@rimath.saitama-u.ac.jp}). Supported in part by JSPS Grant-in-Aid for Scientific Research on Innovative Areas ``Discrete Geometric Analysis for Materials Design" 17H06460.
}
}
\date{}

\begin{document}
\maketitle

\begin{abstract}
In the present article, 
we provide a characterization of a lower $N$-weighted Ricci curvature bound for $N\in ]-\!\infty,1\,]\,\cup\,[\,n,+\infty\,]$ with $\varepsilon$-range introduced by Lu-Minguzzi-Ohta \cite{LMO:CompaFinsler} in terms of a convexity of entropies over Wasserstein space.
We further derive various interpolation inequalities and functional inequalities.
\end{abstract}


 
\section{Introduction}

In this paper,
we present a characterization of a lower $N$-weighted Ricci curvature bound for $N\in ]-\!\infty,1\,]\,\cup\,[\,n,+\infty\,]$ with $\varepsilon$-range introduced by Lu-Minguzzi-Ohta \cite{LMO:CompaFinsler} by a convexity of entropies on the Wasserstein space via mass transport theory.

\subsection{Background}
We first recall the formulation of the weighted Ricci curvature,
and some works on the comparison geometry.
Let $(M,d,\m)$ denote an $n$-dimensional weighted Riemannian
manifold, namely, $M =(M,g)$ is an $n$-dimensional complete Riemannian
manifold, $d$ is 
the Riemannian distance on $M$,
and $\m:=e^{-f}{\vol}_g$ for $f\in C^{\infty}(M)$.
For $N\in]-\infty,\,+\infty\,]$, the associated
{\it $N$-weighted Ricci curvature} $\Ric_f^N$
is defined as follows (\cite{BE1}, \cite{Lich}): 
\begin{align*}
\Ric_f^N:=\Ric_g+{\rm \nabla^2}f-\frac{df\otimes df}{N-n}.
\end{align*}
Here when $N=+\infty$,
we interpret the last term of the right hand side as the limit $0$,
and when $N=n$,
we only consider a constant function $f$, and
set $\Ric_f^n:=\Ric_g$.

It is well-known that
lower weighted Ricci curvature bounds imply various comparison geometric results.
In the classical case of $N\in [\,n,+\infty\,[$,
under a curvature condition
\begin{align}
\Ric_f^N\geq Kg \label{eq:RicciLowerBdd}
\end{align}
for $K\in \mathbb{R}$,
such investigations have been done by \cite{Lo}, \cite{Qi}, \cite{WW},
and so on.

In recent years,
the validity of the $N$-weighted Ricci curvature with $N\in ]-\infty,n[$ has begun to be pointed out (see e.g., \cite{Kl}, \cite{KM}, \cite{KL}, \cite{KwSk}, \cite{KS}, \cite{Lim}, \cite{LMO:CompaFinsler}, \cite{Mai1}, \cite{Mai2}, \cite{Mineg}, \cite{Oh<0}, \cite{Oneedle}, \cite{OT1}, \cite{OT2}, \cite{Sak}, \cite{Sk:OneDim}, \cite{Wy}, \cite{WyYero}).
Wylie-Yeroshkin \cite{WyYero} have proposed a curvature condition
\begin{align}
\Ric_f^1\geq (n-1)\kappa e^{-\frac{4f}{n-1}}g \label{eq:RicciLowerBddWY}
\end{align}
for $\kappa \in \mathbb{R}$ in view of the study of projectively equivalent affine connection,
and established an optimal Laplacian comparison theorem, Bonnet-Myers theorem, Bishop-Gromov volume comparison theorem.
Remark that
before the work of them,
Wylie \cite{Wy} has obtained a splitting theorem of Cheeger-Gromoll type for $\kappa=0$.
For $N\in ]-\infty,1]$,
the first named author and Li \cite{KL} have extended the condition \eqref{eq:RicciLowerBddWY} to
\begin{align}
\Ric_f^N\geq (n-N)\kappa e^{-\frac{4f}{n-N}}g,\label{eq:RicciLowerBddKL}
\end{align}
and generalized the comparison theorems in \cite{WyYero}.

Very recently,
Lu-Minguzzi-Ohta \cite{LMO:CompaFinsler} have suggested a new approach that enables us to investigate the conditions \eqref{eq:RicciLowerBdd} with $K=(N-1)\kappa$, \eqref{eq:RicciLowerBddWY} and \eqref{eq:RicciLowerBddKL} in a unified way. For $N\in ]-\infty,1\,]\,\cup \,[\,n,+\infty\,]$,
they have introduced the notion of the \emph{$\varepsilon$-range}:
\begin{align}
\varepsilon=0\ \text{ for }\ N=1,\qquad \varepsilon\in ]-\sqrt{\varepsilon_0},\sqrt{\varepsilon_0}[  \ \text{ for }\  N\ne1,n,\qquad \varepsilon\in\R \ \text{ for } \ N=n,\label{eq:epsilonrange}
\end{align}
where
\begin{align*}
\varepsilon_0:=\frac{N-1}{N-n}.
\end{align*}
When $N=+\infty$,
we interpret $\eps_0$ as the limit $1$; in particular $\eps \in ]-1,1\,[$.
Within this $\varepsilon$-range,
they have considered a curvature condition
\begin{align}
\Ric_f^N\geq c^{-1}\,\kappa\, e^{-\frac{4(1-\varepsilon)f}{n-1}}g\label{eq:WeightedLowerBddd}
\end{align}
for $\kappa\in\R$, which covers the previous curvature conditions by running $\varepsilon$ over $\varepsilon$-range.
Here $c=c_{N,\varepsilon}\in]0,1]$ is the associated positive constant defined by 
\begin{align}
c:=\frac{1}{n-1}\left(1-\varepsilon^2 \frac{N-n}{N-1} \right)\label{eq:constant}
\end{align}
if $N\ne1$,
and $c:=(n-1)^{-1}$ if $N=1$.
Here we interpret $c$ as the limit $(n-1)^{-1}(1-\varepsilon^2)$ in the case of $N=+\infty$.
When $N\in [\,n,+\infty\,[$ and $\varepsilon=1$ with $c=(N-1)^{-1}$, the curvature condition \eqref{eq:WeightedLowerBddd} covers \eqref{eq:RicciLowerBdd} with $K=(N-1)\kappa$.
Also,
when $N=1$ and $\varepsilon =0$ with $c=(n-1)^{-1}$, it does \eqref{eq:RicciLowerBddWY},
and when $N\in ]-\infty,1]$ and $\varepsilon=\varepsilon_0$ with $c=(n-N)^{-1}$, 
it does \eqref{eq:RicciLowerBddKL}. 
Under the condition \eqref{eq:WeightedLowerBddd},
they have developed comparison geometry in the framework of weighted Finsler manifolds and weighted Finsler space-times.

\subsection{Main results}\label{subsec:MainResults}
Let us introduce our main results.
Lower $N$-weighted Ricci curvature bounds are well-known to be characterized by convexities of entropies on the Wasserstein space.
In the classical case of $N \in [\,n,+\infty\,[$, the characterization of \eqref{eq:RicciLowerBdd} is due to Sturm \cite{StI}, \cite{StII}, and Lott-Villani \cite{LV2}.
Based on such a result, they have independently introduced the so-called 
{\it curvature-dimension condition} ${\rm CD}(K,N)$ for metric measure spaces that is equivalent to \eqref{eq:RicciLowerBdd} in the smooth setting. 
The second named author \cite{Sk:OneDim} gave a characterization of \eqref{eq:RicciLowerBddWY}.

We now aim to provide a characterization of the curvature condition \eqref{eq:WeightedLowerBddd}.
Let $N\in ]-\infty,1\,]\,\cup \,[\,n,+\infty\,]$, 
and $\varepsilon \in \mathbb{R}$ in the range \eqref{eq:epsilonrange}.
Let $\mathcal{DC}_{N,\varepsilon}$ be the set of all continuous convex functions $U : [\,0,+\infty\,[\to\R$ with $U(0) = 0$
such that a function $\varphi_U :\,]\,0,+\infty\,[\to\R$ defined by $\varphi_U(r) := r^{\frac{c+1}{c}} U(r^{-\frac{c+1}{c}})$ is
convex,
where the constant $c>0$ is defined as \eqref{eq:constant} ($\mathcal{DC}_{N,\varepsilon}$ can be written as $\mathcal{DC}_{(c+1)/c}$ in the notation of \cite{Vi2}).
Let $\mathscr{P}_2(M)$ be the set of all Borel probability measures on $M$ with finite second moment,
which is endowed with the \textit{$L^{2}$-Wasserstein distance function} $W_2$.
For $U\in\mathcal{DC}_{N,\varepsilon}$, a functional $U_{\m}$ on $\mathscr{P}_2(M)$ is defined by \begin{align}
U_{\m}(\mu):=\int_MU(\rho)\d\m,\label{eq:Displacement}
\end{align}
where $\rho$ is the density of the absolutely continuous part in the Lebesgue
decomposition of $\mu$ with respect to $\m$. For a function $H\in\mathcal{DC}_{N,\varepsilon}$ defined 
by $H(r) := c^{-1}(c+1)r(1-r^{-\frac{c}{c+1}})$, the functional $H_{\m}$ on $\mathscr{P}_2(M)$ defined as \eqref{eq:Displacement}
is called the {\it R\'enyi entropy}.

Following \cite{Sk:OneDim},
we introduce a twisted coefficient in our setting.
We define two lower semi continuous functions $d_{N,\varepsilon,f,t}, d_{N,\varepsilon,f} : M \times M\to\R$ by 
\begin{align*}
d_{N,\varepsilon,f,t}(x,y):=\inf_{\gamma}\int_0^{td(x,y)}e^{-\frac{2(1-\varepsilon)f(\gamma(\xi))}{n-1}}\d \xi,\qquad d_{N,\varepsilon,f}:=d_{N,\varepsilon,f,1}
\end{align*}
for $t\in[0,1]$,
where the infimum is taken over all unit speed minimal geodesics $\gamma:[0,d(x,y)]\to M$ from $x$ to $y$. The function $d_{N,\varepsilon,f}$ is called the 
{\it re-parametrized distance} (cf. \cite{WyYero}).
Note that for $t\in]\,0,1\,[$, the function $d_{N,\varepsilon,f,t}$ 
is not always symmetric. For  $\kappa\in\R$, let 
$\s_{\kappa}(s)$ stand for a unique solution of
the Jacobi equation  $\psi''(s) + \kappa \psi(s) = 0$ 
with $\psi(0) = 0, \psi'(0) = 1$, and $C_{\kappa}$ the diameter of the space form of constant curvature $\kappa$.
More precisely,
they can be written as
\begin{equation*}
\mathfrak{s}_{\kappa}(s)=\begin{cases}
                                                     \displaystyle \frac{\sin \sqrt{\kappa}s}{\sqrt{\kappa}} & \text{if $\kappa>0$}, \\
                                                                                s          & \text{if $\kappa=0$},\\
                                                     \displaystyle \frac{\sinh \sqrt{\vert\kappa \vert}s}{\sqrt{\vert \kappa \vert}}           & \text{if $\kappa<0$},
                                                   \end{cases}\quad 
C_{\kappa}=\begin{cases}
                                                     \displaystyle \frac{\pi}{\sqrt{\kappa}} & \text{if $\kappa>0$}, \\
                                                                                \infty           & \text{if $\kappa \leq 0$}.
                                                   \end{cases}
\end{equation*}
For $t\in ]\,0,1\,[$, we define the {\it twisted coefficient $\beta_{\kappa,N,\varepsilon,f,t}:M\times M \to \mathbb{R}\cup \{+\infty\}$} by 
\begin{align}
\beta_{\kappa,N,\varepsilon,f,t}(x,y):=\left(\frac{\;\s_{\kappa}(d_{N,\varepsilon,f,t}(x,y))\;}{t \,\s_{\kappa}(d_{N,\varepsilon,f}(x,y))} \right)^{c^{-1}} \label{eq:twistedcoeffi}
\end{align}
if $d_{N,\varepsilon,f}(x,y)\in ]0,C_{\kappa}[$;  $\beta_{\kappa,N,\varepsilon,f,t}(x,y)=1$ if $x=y$; otherwise, $\beta_{\kappa,N,\varepsilon,f,t}(x, y) := +\infty$. 

\begin{remark}\label{rem:TwistedCoeffic}
{\rm  The definition of the twisted coefficient for $x=y$ is reasonable since we see $\beta_{\kappa,N,\varepsilon,f,t}(x,y)\to 1$ as $d(x,y)\to 0$ (see Appendix for the proof).}
\end{remark}
Let $\mathscr{P}_2^{\rm ac}(M)$ denote the set of all Borel probability measures in $\mathscr{P}_2(M)$ 
that are absolutely 
continuous with respect to $\m$.
We now introduce the following convexity properties:

\begin{defn}\label{defn:LV}
{\rm  Let $\kappa\in\R,\,N\in ]-\infty,1\,]\,\cup \,[\,n,+\infty\,]$,
and $\varepsilon \in \mathbb{R}$ in the range \eqref{eq:epsilonrange}.
We say that $(M,d,\m)$ satisfies {\it the twisted curvature-dimension condition ${\rm TwCD}(\kappa,N,\varepsilon)$} if for every pair $\mu_0,\mu_1\in \mathscr{P}_2^{\rm ac}(M)$,
\begin{align}
U_{\m}(\mu_t)&\leq (1-t)\int_{M^2}U\left(\frac{\rho_0(x)}{{\beta}_{\kappa,N,\varepsilon,f,1-t}(y,x)} \right)\frac{\;{\beta}_{\kappa,N,\varepsilon,f,1-t}(y,x)\;}{\rho_0(x)}\pi(\d x\d y)\label{eq:kapaptwisted}\\
&\hspace{3cm}+t\int_{M^2}U\left(\frac{\rho_1(y)}{{\beta}_{\kappa,N,\varepsilon,f,t}(x,y)} \right)\frac{\;{\beta}_{\kappa,N,\varepsilon,f,t}(x,y)\;}{\rho_1(y)}\pi(\d x\d y)\nonumber
\end{align}
for all $U\in\mathcal{DC}_{N,\varepsilon}$ and $t\in]\,0,1\,[$,
where $\rho_i$ is the density of $\mu_i$ with respect to $\m$ for each $i=0,1$,
and $\pi$ is a unique optimal coupling of $(\mu_0,\mu_1)$,
and $(\mu_t)_{t\in[0,1]}$ is a unique minimal geodesic in the $L^2$-Wasserstein space $(\mathscr{P}_2(M),W_2)$ from $\mu_0$ to $\mu_1$, which lies in $\mathscr{P}_2^{\rm ac}(M)$.
}
\end{defn}

\begin{remark}\label{rem:compact}
{\rm In Definition \ref{defn:LV}, we only consider $U \in \mathcal{DC}_{N,\varepsilon}$ such that \eqref{eq:kapaptwisted} makes sense for all $\mu_0,\mu_1\in \mathscr{P}_2^{\rm ac}(M)$.
Notice that
for $H\in\mathcal{DC}_{N,\varepsilon}$ defined as $H(r):=c^{-1}(c+1)r(1-r^{-\frac{c}{c+1}})$,
\eqref{eq:kapaptwisted} makes sense for all $\mu_0,\mu_1$.
For general $U \in \mathcal{DC}_{N,\varepsilon}$,
such a property is guaranteed by a condition
\begin{equation*}
\int_{M}\,\frac{1}{(1+d(x,x_0)^2)^{1/c}}\,\m(\d x)<+\infty
\end{equation*}
for some $x_0\in M$ concerning the reference measure $\m$ (cf. \cite[Theorems 17.8, 17.28]{Vi2}).}
\end{remark}

\begin{defn}
{\rm  Let $\kappa\in\R,\,N\in ]-\infty,1\,]\,\cup \,[\,n,+\infty\,]$,
and $\varepsilon \in \mathbb{R}$ in the range \eqref{eq:epsilonrange}. We say that $(M,d,\m)$ satisfies
{\it the relaxed twisted curvature-dimension condition ${\rm TwCD_{rel}}(\kappa,N,\varepsilon)$} if the inequality \eqref{eq:kapaptwisted} holds for $H\in\mathcal{DC}_{N,\varepsilon}$ 
defined as $H(r):=c^{-1}(c+1)r(1-r^{-\frac{c}{c+1}})$.
}
\end{defn}


\begin{remark}\label{rem:LVSt}
{\rm In the case of $N\in [\,n,+\infty\,[$ and $\varepsilon=1$, the condition ${\rm TwCD}(\kappa,N,1)$ coincides with the curvature-dimension
condition ${\rm CD}((N-1)\kappa,N)$ in the sense of Lott-Villani \cite{LV2}.
Similarly, ${\rm TwCD_{rel}}(\kappa,N,1)$ coincides
with ${\rm CD}((N-1)\kappa,N)$ in the
sense of Sturm \cite{StI}, \cite{StII}.
In the case of $N=1$ with $\varepsilon=0$, the conditions ${\rm TwCD}(\kappa,1,0)$ and ${\rm TwCD_{rel}}(\kappa,1,0)$ coincide with the $\kappa$-twisted curvature bound and the relaxed one in \cite{Sk:OneDim},
respectively.}
\end{remark}

We now state our main theorem.

\begin{thm}\label{thm:Main}
Let $\kappa\in\R,\,N\in ]-\infty,1\,]\,\cup\, [\,n,+\infty\,]$,
and $\varepsilon \in \mathbb{R}$ in the range \eqref{eq:epsilonrange}.
We additionally assume that if $N\neq 1,n$, then $\eps \neq 0$.
Then the following are equivalent:
\begin{enumerate}
\item\label{item:Main1} $\Ric_f^N\geq c^{-1}\,\kappa\, e^{-\frac{4(1-\varepsilon)f}{n-1}}g$;
\item\label{item:Main2} $(M,d,\m)$ satisfies ${\rm TwCD}(\kappa,N,\varepsilon)$;
\item\label{item:Main3} $(M,d,\m)$ satisfies ${\rm TwCD_{rel}}(\kappa,N,\varepsilon)$.
\end{enumerate}
\end{thm}

\begin{remark}
The restriction $\eps \neq 0$ is due to a technical issue,
which is not a natural requirement.
The authors do not know whether this can be removed.
\end{remark}

In the case of $N\in [\,n,+\infty\,[$ and $\varepsilon=1$,
Theorem \ref{thm:Main} is nothing but the well-known characterization of the curvature condition (\ref{eq:RicciLowerBdd}) with $K=(N-1)\kappa$ by ${\rm CD}((N-1)\kappa,N)$ (see \cite[Theorem 4.22]{LV2}, and also \cite[Theorem 1.7]{StII}).
When $N\in [\,n,+\infty\,[$,
Theorem \ref{thm:Main} for $\varepsilon \neq 1$ is new and not treated in the literature.

The second named author \cite{Sk:OneDim} has shown Theorem \ref{thm:Main} when $N=1$ (see \cite[Theorem~1.4]{Sk:OneDim}).
Theorem \ref{thm:Main} for $N\in ]-\infty,1\,[$ is a new result;
in particular,
by letting $\varepsilon=\varepsilon_0$,
one can obtain the following characterization of the condition (\ref{eq:RicciLowerBddKL}):
\begin{cor}\label{cor:Main2}
Under the same setting as in Theorem \ref{thm:Main},
if $N\in ]-\infty,1]$,
then the following statements are equivalent: 
\begin{enumerate}
\item\label{item:Main1} $\Ric_f^N\geq (n-N)\,\kappa\, e^{-\frac{4f}{n-N}}g$;
\item\label{item:Main2} $(M,d,\m)$ satisfies ${\rm TwCD}(\kappa,N,\varepsilon_0)$;
\item\label{item:Main3} $(M,d,\m)$ satisfies ${\rm TwCD_{rel}}(\kappa,N,\varepsilon_0)$.
\end{enumerate}
\end{cor}

We also notice that
as a corollary of the proof of Theorem \ref{thm:Main},
we obtain the following (see Proposition \ref{prop:CDimpliesConvexity} below):
\begin{cor}\label{cor:Main}
Under the same setting as in Theorem \ref{thm:Main},
the implication from \ref{item:Main1} to \ref{item:Main2} always holds $($without the restriction $\varepsilon \neq 0$$)$. 
\end{cor}

From this viewpoint,
under the curvature condition \eqref{eq:WeightedLowerBddd},
we derive several interpolation inequalities such as
$p$-mean inequality,
Pr\'ekopa-Leindler inequality,
Borell-Brascamp-Lieb inequality,
and Brunn-Minkowski inequality (see Subsection \ref{subsec:Inteplation}),
and also study functional inequalities (see Section \ref{sec:FuncIneq}). 

\section{Preliminaries}\label{sec:Preliminary}
This section is devoted to basics on optimal transport theory and comparison geometry.

\subsection{Optimal transport theory}\label{subsec:OptomalTransport}
We recall some basic facts on the optimal transport theory.
Referring to \cite{CMS1}, \cite{Oint}, \cite{Vi2},
we use the same notation and terminology as in the preliminaries of \cite{Sk:OneDim} (see \cite[Subsection~2.2]{Sk:OneDim}).
On a metric space $(Z,d_Z)$,
a curve $\gamma : [0,l] \to  Z$ is said to be a {\it minimal geodesic} if there is $a \geq 0$ such that $d_Z(\gamma(t_0),\gamma(t_1)) = a|t_0 - t_1|$ for all $t_0,t_1 \in[0,l]$. Moreover, if $a = 1$, then $\gamma$ is said to be a {\it unit speed minimal geodesic}.

Let $\mathscr{P}(M)$ be the set of all Borel probability measures on $M$.
For $\mu,\nu\in\mathscr{P}(M)$, a Borel probability measure $\pi$ on $M\times M$  is said to be a {\it coupling of $(\mu,\nu)$} if $\pi(X\times M)=\mu(X)$ 
and $\pi(M \times X) = \nu(X)$ for all Borel subsets 
$X \subset  M$. Let $\Pi(\mu,\nu)$ stand for the set of all 
couplings of $(\mu,\nu)$.
Recall that $\mathscr{P}_2(M)$ denotes the set of all Borel probability measures on $M$ with finite second moment,
namely,
$\mu\in \mathscr{P}_2(M)$ if
$$
\int_Md(x,x_0)^2\mu(\d x) < +\infty
$$ for 
some $x_0\in M$.
The \emph{$L^2$-Wasserstein distance function} $W_2$ is defined as  
\begin{align}
W_2(\mu,\nu):=\inf_{\pi\in\Pi(\mu,\nu)}\left(\int_{M^2}d(x,y)^2\pi(\d x\d y) \right)^{\frac12}.\label{eq:Wasserstein}
\end{align}
The pair $(\mathscr{P}_2(M),W_2)$ is known to be a complete separable metric space (see e.g., \cite[Theorem 6.18]{Vi2}), and called the {\it $L^2$-Wasserstein space}. 
A coupling $\pi\in\Pi(\mu,\nu)$ is said to be \emph{optimal} if it attains the infimum of \eqref{eq:Wasserstein}.
Recall the following fundamental result on the optimal coupling in smooth setting due to Brenier \cite{Brenier:Polar1991}, McCann \cite{McCann:Polar}, and Figalli-Gigli \cite{FigalliGigli:2011}  (see \cite{Brenier:Polar1991}, \cite[Theorem 1]{FigalliGigli:2011}, \cite[Theorem 3]{McCann:Polar}):
\begin{thm}\label{thm:BrenierMcCann}
For $\mu\in\mathscr{P}^{
ac}_{2}(M)$ and  $\nu\in\mathscr{P}_2(M)$, there is a locally semi-convex function 
$\phi$ on an open subset $\Omega$ of $M$ with $\mu(\Omega)=1$ such that a map $F_t$ defined by
\begin{align}
F_t(z) := \exp_z(t\nabla\phi(z))\label{eq;OptimalMap}
\end{align}
provides a unique optimal coupling $\pi$ of $(\mu,\nu)$ via the pushforward measure $\pi:=(F_0\times F_1)_{\sharp}\mu$ of $\mu$ by $F_0\times F_1$, and also determines a unique minimal geodesic $(\mu_t)_{t\in[0,1]}$ in $(\mathscr{P}_2(M),W_2)$ from $\mu$ to $\nu$ via $\mu_t:=(F_t)_{\sharp}\mu_0$.
\end{thm}

The function $\phi$ provided in 
Theorem~\ref{thm:BrenierMcCann} is called the {\it Kantorovich potential}, which is twice differentiable $\mu$-almost everywhere as a consequence of the Alexandrov-Bangert theorem.
The Kantorovich potential $\phi$ has the following properties (see \cite[Theorem 1.1]{FigalliGigli:2011}, \cite[Proposition 4.1, Corollary 5.2]{CMS1}): If $\phi$ is twice differentiable at $x$, then $F_t(x)$ does not belong to the cut locus ${\rm Cut}(x)$ of $x$, and the differential
$(dF_t)_x$ is well-defined for every
$t\in[0,1]$. Also, $\phi$ satisfies the following (see \cite[Theorem 8.7]{Vi2}): The curve $(\mu_t)_{t\in[0,1)}$ lies in $\mathscr{P}_2^{\rm ac}(M)$. 
We finally recall the \emph{Monge-Amp\`{e}re equation} (see \cite[Theorem 11.1]{Vi2}):

\begin{thm}\label{thm:CorderoMcCannSche}
Let $\mu,\nu\in\mathscr{P}_2^{ac}(M)$, 
and let $\phi$ be the Kantorovich 
potential obtained in Theorem~\ref{thm:BrenierMcCann}. Then for $\mu$-almost every $x$, we have:
\begin{enumerate}
\item $\phi$ is twice differentiable at $x$;
\item the determinant ${\rm det}(dF_t)_x$ is positive for every $t\in[0,1]$;
\item $\rho_0(x)=\rho_1(F_1(x))e^{-f(F_1(x))+f(x)}{\rm det}(dF_1)_x$, where $\rho_0$ and $\rho_1$ are the densities of $\mu$ and of $\nu$ with respect to $\m$, respectively.
\end{enumerate} 
\end{thm}

\subsection{Comparison geometric results}
We next review one of comparison geometric results,
which will be used in the proof of the main theorem.
Let $N\in ]-\infty,1\,]\,\cup\, [\,n,+\infty\,]$, 
and $\varepsilon \in \mathbb{R}$ in the range \eqref{eq:epsilonrange}.

For $x\in M$,
let $U_{x}M$ be the unit tangent sphere at $x$.
For $v\in U_{x}M$,
let $\gamma_{v}:[0,\infty)\to M$ denote the unit speed geodesic with initial conditions $\gamma_{v}(0)=x$ and $\dot{\gamma}_{v}(0)=v$.
Define a function $s_{N,\varepsilon,f,v}:[0,+\infty\,]\to [0,s_{N,\varepsilon,f,v}(+\infty)]$ by
\begin{equation*}
s_{N,\varepsilon,f,v}(t):=\int^{t}_{0}\,e^{-\frac{2(1-\varepsilon)f(\gamma_{v}(\xi))}{n-1}}\, \d \xi.
\end{equation*}
We also set
\begin{align}
\tau(v):=\sup\{t>0\mid d(x,\gamma_v(t))=t\}, \quad \tau_{N,\varepsilon,f}(v):=s_{N,\varepsilon,f,v}(\tau(v)).\label{eq:directionalDiam}
\end{align}
The authors \cite{KwSk} has shown the following (see \cite[Lemma~2.6, Proposition~3.1]{KwSk}):
\begin{prop}\label{prop:DirectionalMyers}
For $\kappa>0$,
if \,$\Ric_f^N\geq c^{-1}\,\kappa \,e^{-\frac{4(1-\varepsilon)f}{n-1}}g$, then for
all $x\in M$ and $v\in U_xM$, 
\begin{align*}
\tau_{N,\varepsilon,f}(v) \leq C_{\kappa}.
\end{align*}
Moreover, for the re-parametrized distance $d_{N,\varepsilon,f}$,
we have
\begin{align*}
\sup_{x,y\in M}d_{N,\varepsilon,f}(x,y)\leq C_{\kappa}.
\end{align*}
\end{prop}

Note that
the authors \cite{KwSk} have obtained a similar comparison result in a more general setting such that the density is a vector field, and $\kappa$ is variable.

\section{Key inequalities}\label{sec:keyIneq}
Hereafter,
we always fix $N\in ]-\infty,1\,]\,\cup\, [\,n,+\infty\,]$,
and $\varepsilon \in \mathbb{R}$ in the range \eqref{eq:epsilonrange}.
Moreover,
in the case of $N=n$,
the density function $f$ is constant;
in particular,
the main assertions have been already proved in the works of Sturm \cite{StI}, \cite{StII}, and Lott-Villani \cite{LV2}.
Furthermore,
in the case of $N=1$,
they have been done by the second named author \cite{Sk:OneDim}.
Thus, we further suppose $N\neq 1,n$. 

The aim of this section is to produce the following key inequality for the proof of our main theorem (cf. \cite[Proposition 3.1]{Sk:OneDim}):

\begin{prop}\label{prop:Jacobian}
Let $\mu,\nu\in \mathscr{P}_2^{\rm ac}(M)$, and let $\phi$ be the Kantorovich potential in Theorem~\ref{thm:BrenierMcCann}. 
For a fixed $x\in M$,
assume that $\phi$ is twice differentiable at $x$, and ${\rm det}\,(dF_t)_x>0$ for every $t\in[0,1]$,
where $F_t$ is defined as \eqref{eq;OptimalMap}.
For each $t\in[0,1]$, set
\begin{align}
J_t(x):=e^{-f(F_t(x))+f(x)}{\rm det}(dF_t)_x.\label{eq:Jacobian}
\end{align}
For $\kappa\in\R$, if ${\Ric}_f^N\geq c^{-1}\,\kappa\, e^{-\frac{4(1-\varepsilon)f}{n-1}}g$, then for every $t\in]\,0,1\,[$
\begin{align*}
J_t(x)^{\frac{c}{c+1}}\geq (1-t){\beta}_{\kappa,N,\varepsilon,f,1-t}(F_1(x),x)^{\frac{c}{c+1}}J_0(x)^{\frac{c}{c+1}}+
t{\beta}_{\kappa,N,\varepsilon,f,1-t}(x,F_1(x))^{\frac{c}{c+1}}{J_1(x)}^{\frac{c}{c+1}}.
\end{align*} 
\end{prop}

Throughout this section, let $\mu,\nu,\phi,x$ be as in Proposition~\ref{prop:Jacobian}.

\subsection{Riccati inequalities}\label{subsec:RiccatiEq}
Define a curve $\gamma:[0,1] \to M$ by $\gamma(t):=F_{t}(x)$,
and choose an orthonormal basis $\{e_{i}\}^{n}_{i=1}$ at $x$ with $e_{n}=\dot{\gamma}(0)/\Vert \dot{\gamma}(0)\Vert$.
For each $i$,
we define a Jacobi field $E_{i}$ along $\gamma$ by $E_{i}(t):=(dF_{t})_{x}(e_{i})$.
For each $t \in [0,1]$
let $A(t)=(a_{ij}(t))$ be an $n \times n$ matrix determined by
\begin{equation*}
E'_{i}(t)=\sum^{n}_{j=1}\,a_{ij}(t)\,E_{j}(t).
\end{equation*}
Let us consider a function $h:[0,1]\to \mathbb{R}$ defined by
\begin{equation*}
h(t):=\log \,\det (dF_{t})_{x} -\int^{t}_{0}\,a_{nn}(\xi)\,\d\xi,
\end{equation*}
which enjoys the following Riccati inequality (see e.g., (1.4), (1.9) in \cite{StII}, and (14.21) in \cite{Vi2}):
 
\begin{lem}\label{lem:RiccatiIneq}
For every $t\in]\,0,1\,[$ we have 
\begin{align*}
h''(t)\leq-\frac{h'(t)^2}{n-1}-\Ric_g(\dot{\gamma}(t)). 
\end{align*}
\end{lem}

We define a function $l : [0,1]\to\R$ by
\begin{align*}
l(t):=h(t)-f(\gamma(t))+f(x).
\end{align*}   

We show the following Riccati inequality,
which is compatible with our setting (cf. \cite[Lemma 3.3]{Sk:OneDim}, and also \cite[Lemma~2.1]{KwSk} in the literature of comparison geometry). 
\begin{lem}\label{lem:JacocbianIneq0}
For every $t\in]\,0,1\,[$ we have 
\begin{align}
\left(e^{\frac{2(1-\varepsilon)f(\gamma(t))}{n-1}}l'(t)\right)'\leq 
-e^{\frac{2(1-\varepsilon)f(\gamma(t))}{n-1}}
\left(c\,l'(t)^2+\Ric_f^N(\dot{\gamma}(t) \right).\label{eq:GeneBishop}
\end{align}
\end{lem}
\begin{pf}
Set $f_{x}:=f\circ \gamma$.
Lemma \ref{lem:RiccatiIneq} leads us to
\begin{align*}
l''(t)=h''(t)-f''_{x}(t) & \leq -\frac{h'(t)^{2}}{n-1} -\left(\Ric_{g}(\dot{\gamma}(t))+ f''_{x}(t)   \right)\\
                            &   =   -c\,l'(t)^2-\frac{2(1-\varepsilon)l'(t)\,f'_{x}(t)}{n-1}-\Ric^{N}_{f}(\dot{\gamma}(t))\\
                            &\hspace{1cm}-\frac{1}{n-1}\left(\varepsilon \sqrt{\frac{N-n}{N-1}}l'(t)+ \sqrt{\frac{N-1}{N-n}}f'_{x}(t)   \right)^2\\
                            &   \leq   -c\,l'(t)^2-\frac{2(1-\varepsilon)l'(t)\,f'_{x}(t)}{n-1}-\Ric^{N}_{f}(\dot{\gamma}(t)).
\end{align*}
This implies
\begin{align*}
e^{\frac{-2(1-\varepsilon)f_{x}(t)}{n-1}}\left(   e^{\frac{2(1-\varepsilon)f_{x}(t)}{n-1}} \,l'(t)    \right)'      =   l''(t)+\frac{2(1-\varepsilon)\,l'(t)\,f'_{x}(t)}{n-1}\leq -c\,l'(t)^{2} -\Ric^{N}_{f}(\dot{\gamma}(t)).
\end{align*}
We arrive at the desired inequality \eqref{eq:GeneBishop}.
\end{pf}

\subsection{Jacobian inequalities}\label{subsec:Jacobian}

Once we obtain the Riccati inequality \eqref{eq:GeneBishop},
one can prove Proposition \ref{prop:Jacobian} by the same argument as in the proof of \cite[Proposition 3.1]{Sk:OneDim}.
Define a function $D:[0,1]\to\R$ by 
\begin{align*}
D(t):=\exp\left(c\,l(t) \right).
\end{align*}
In virtue of Lemma~\ref{lem:JacocbianIneq0},
we have the following (cf. \cite[Lemma 3.5]{Sk:OneDim}):
\begin{lem}\label{lem:JacocbianIneq2}
For  $\kappa\in\R$, if ${\Ric}_f^N\geq c^{-1}\kappa\, e^{-\frac{4(1-\varepsilon)f}{n-1}}g$, then for every $t\in]\,0,1\,[$ we have 
\begin{align*}
D(t)\geq \frac{\s_{\kappa}(d_{N,\varepsilon,f,1-t}(F_1(x),x))}{\s_{\kappa}(d_{N,\varepsilon,f}(F_1(x),x))}D(0)+\frac{\s_{\kappa}(d_{N,\varepsilon,f,t}(x,F_1(x)))}{\s_{\kappa}(d_{N,\varepsilon,f}(x,F_1(x)))}D(1).
\end{align*}
\end{lem}
\begin{pf}
As in the proof of \cite[Lemma 3.5]{Sk:OneDim},
we define $s_{f}:[0,1]\to \mathbb{R}$ by
\begin{align*}
s_{f}(t):=\int_0^te^{-\frac{2(1-\varepsilon)f(\gamma(\xi))}{n-1}}\d \xi.
\end{align*}
For $a:=s_{f}(1)$,
we further define $\wh{l},\wh{D}:[0,a]\to\R$ by
\begin{align*}
\wh{l}:=l \circ t_{f},\quad \wh{D}:=D\circ t_{f},
\end{align*}
where $t_{f}:[0,a]\to[0,1]$ is the inverse function of $s_{f}$.
For each $s\in]0,a[$ it holds that
\begin{align}
c^{-1}\frac{\wh{D}''(s)}{\wh{D}(s)}=\wh{l}''(s)+c\,\wh{l}'(s)^2.\label{eq:RiccatiIneq1}
\end{align}
We also define functions $L:[0,1]\to\R$ and $\wh{L}:[0,a]\to\R$ by
\begin{align*}
L(t):=e^{\frac{2(1-\varepsilon)f(\gamma(t))}{n-1}}l'(t),\quad \wh{L}:=L\circ t_{f}.
\end{align*}
By Lemma~\ref{lem:JacocbianIneq0} we see
\begin{align}\label{eq:RiccatiIneq2}
\widehat{l}''(s)   &    =    \widehat{L}'(s)=t'_{f}(s)\,L'(t_{f}(s))\\ \notag
                        &   \leq  -e^{\frac{4(1-\varepsilon)f\left(\gamma\left(t_{f}(s)\right)\right)}{n-1}}\left(      c\, l'(t_{f}(s))^{2}+\Ric^{N}_{f}(\dot{\gamma}(t_{f}(s))) \right)\\  \notag
                        &     =   -c\,\widehat{l}'(s)^{2}-e^{\frac{4(1-\varepsilon)f\left(\gamma\left(t_{f}(s)\right)   \right)}{n-1}}  \,\Ric^{N}_{f}(\dot{\gamma}(t_{f}(s))). \notag
\end{align}
The equality \eqref{eq:RiccatiIneq1} together with \eqref{eq:RiccatiIneq2} yields
\begin{align*}
c^{-1}\frac{\wh{D}''(s)}{\wh{D}(s)}\leq 
-e^{\frac{4(1-\varepsilon)f\left(\gamma(t_{f}(s))\right)}{n-1}}
\Ric^{N}_{f}(\dot{\gamma}(t_{f}(s)))\leq -c^{-1}\,\kappa\, d(x,y)^2,
\end{align*}
where $y:=F_1(x)$. Hence, $\wh{D}''(s)+\kappa\, d(x,y)^2\wh{D}(s)\leq0$ on $]0,a[$. 

Since $\phi$ is twice differentiable at $x$, the curve $\gamma$ lies in the complement of ${\rm Cut}(x)$. In particular, $\gamma$ is a unique minimal geodesic from $x$ to $y$.
Therefore,
\begin{align*}
a\,d(x,y)=d_{N,\varepsilon,f}(x,y)<\tau_{N,\varepsilon,f}\left(\frac{\dot{\gamma}(0)}{\|\dot{\gamma}(0)\|} \right),
\end{align*}
where $\tau_{N,\varepsilon,f}$ is defined as \eqref{eq:directionalDiam}. 
Due to Proposition~\ref{prop:DirectionalMyers}, $\kappa\, d(x,y)^2\in]-\infty,a^{-2}\pi^2[$.
Now,
an elementary comparison argument implies the following (see e.g., \cite[Theorem 14.28]{Vi2}):
For all $s_0,s_1\in[0,a]$ and $\lambda\in[0,1]$, 
\begin{align*}
\wh{D}((1-\lambda)s_0+\lambda s_1)\geq 
\frac{\s_{\kappa}((1-\lambda)|s_0-s_1|d(x,y))}{\s_{\kappa}(|s_0-s_1|d(x,y))}\wh{D}(s_0)
+ \frac{\s_{\kappa}(\lambda|s_0-s_1|d(x,y))}{\s_{\kappa}(|s_0-s_1|d(x,y))}\wh{D}(s_1).
\end{align*}
This implies that
for every $s\in\,]\,0,a\,[$
we also see
\begin{align*}
\wh{D}(s)\geq \frac{\s_{\kappa}((a-s)\,d(x,y))}{\s_{\kappa}(a\,d(x,y))}\wh{D}(0)+\frac{\s_{\kappa}(s\,d(x,y))}{\s_{\kappa}(a\,d(x,y))}\wh{D}(a).
\end{align*}
It follows that
for every $t\in\,]\,0,1\,[$ 
\begin{align*}
D(t)\geq \frac{\s_{\kappa}((a-s_{f}(t))\,d(x,y))}{\s_{\kappa}(a\,d(x,y))}D(0)+
\frac{\s_{\kappa}(s_{f}(t)\,d(x,y))}{\s_{\kappa}(a\,d(x,y))}D(1).
\end{align*}
In view of the uniqueness of the geodesic $\gamma$, for every 
$t\in [0,1]$ it holds that
\begin{align*}
(a-s_{f}(t))\,d(x,y)=d_{N,\varepsilon,f,1-t}(y,x),\quad s_{f}(t)\,d(x,y)=d_{N,\varepsilon,f,t}(x,y).
\end{align*}
Thus,
we complete the proof.
\end{pf}

Let us give a proof of Proposition~\ref{prop:Jacobian}.

\begin{apf}{Proposition~\ref{prop:Jacobian}}
For $\kappa\in\R$, we assume ${\Ric}_f^N\geq c^{-1}\kappa\, e^{-\frac{4(1-\varepsilon)f}{n-1}}g$. 
Set
\begin{align*}
\overline{D}(t):=\exp\left(\int_0^ta_{nn}(\xi)\d \xi \right).
\end{align*}
The following is well-known (see e.g., (1.10) in \cite{StII}, and (14.19) in \cite{Vi2}):
\begin{align}
\overline{D}(t)\geq (1-t)\overline{D}(0)+t\overline{D}(1).
\label{eq:Concavity}
\end{align}
From Lemma~\ref{lem:JacocbianIneq2}, \eqref{eq:Concavity}  
and the H\"older inequality, it follows that
\begin{align*}
J_t(x)^{\frac{c}{c+1}}&=
D(t)^{1-\frac{c}{c+1}}
\overline{D}(t)^{\frac{c}{c+1}}\\
&\geq(1-t)^{\frac{c}{c+1}}
\left(\frac{\s_{\kappa}(d_{N,\varepsilon,f,1-t}(F_1(x),x))}{\s_{\kappa}(d_{N,\varepsilon,f}(F_1(x),x))} \right)^{1-\frac{c}{c+1}}J_0(x)^{\frac{c}{c+1}}\\
&\hspace{2cm}+t^{\frac{c}{c+1}}
\left(\frac{\s_{\kappa}(d_{N,\varepsilon,f,t}(F_1(x),x))}{\s_{\kappa}(d_{N,\varepsilon,f}(F_1(x),x))} \right)^{1-\frac{c}{c+1}}J_1(x)^{\frac{c}{c+1}}.
\end{align*}
This proves Proposition~\ref{prop:Jacobian}.
\end{apf}

\section{Displacement convexity}\label{sec:DisplacementConvex}
In this section, we prove Theorem~\ref{thm:Main} with the help of Proposition~\ref{prop:Jacobian}.

\subsection{Curvature bounds imply displacement convexity}\label{subsec:CDDisplacement}
We first show the implication from \ref{item:Main1} to \ref{item:Main2} in Theorem~\ref{thm:Main},
which is also stated as Corollary \ref{cor:Main} in Subsection \ref{subsec:MainResults} (cf. \cite[Proposition 4.1]{Sk:OneDim} for the case of $N=1$).

\begin{prop}\label{prop:CDimpliesConvexity}
For $\kappa\in\R$, if $\Ric_f^N\geq c^{-1}\kappa\,e^{-\frac{4(1-\varepsilon)f}{n-1}}g$ holds, then $(M,d,\m)$ satisfies ${\rm TwCD}(\kappa,N,\varepsilon)$.
\end{prop}
\begin{pf}
Let $\mu,\nu \in\mathscr{P}_2^{\rm ac}(M)$,
and let $\phi$ be the Kantorovich potential obtained in Theorem~\ref{thm:BrenierMcCann}.
The map $F_t$ defined as \eqref{eq;OptimalMap} provides a unique optimal coupling $\pi$ of $(\mu,\nu)$ via
$\pi := (F_0\times F_1)_{\sharp}\mu$. It also determines a unique minimal geodesic $(\mu_t)_{t\in[0,1]}$ from $\mu$ to $\nu$ via $\mu_t 
:= (F_t)_{\sharp}\mu$. Moreover, 
thanks to Theorem~\ref{thm:CorderoMcCannSche}, for a fixed $t\in]\,0,1\,[$, the Monge-Amp\`ere equations 
\begin{align}
\rho_0(x)=\rho_1(F_1(x))J_1(x)=\rho_t(F_t(x))J_t(x)\label{eq:JacobianEq}
\end{align}
hold for $\mu_0$-almost every $x\in M$, where $\rho_t$ denotes the density of $\mu_t$ with respect to $\m$.
For $U\in\mathcal{DC}_{N,\varepsilon}$, let $\varphi_U(r):=r^{\frac{c+1}{c}}U(r^{-\frac{c+1}{c}})$. 
From \eqref{eq:JacobianEq} and Proposition~\ref{prop:Jacobian},
we deduce
\begin{align*}
U_{\m}(\mu_t)&=\int_MU\left(\frac{\rho_0(x)}{J_t(x)} \right)\frac{J_t(x)}{\rho_0(x)}\mu_0(\d x)
=\int_M\varphi_U\left(\left(\frac{J_t(x)}{\rho_0(x)} \right)^{\frac{c}{c+1}} \right)\mu_0(\d x)
\\
&\leq (1-t)\int_M\varphi_U\left(
{\beta}_{\kappa,N,\varepsilon,f,1-t}(F_1(x),x)^{\frac{c}{c+1}} 
\left(\frac{J_0(x)}{\rho_0(x)} \right)^{\frac{c}{c+1}}
\right)\mu_0(\d x)\\
&\hspace{2cm}
+t\int_M\varphi_U\left(
{\beta}_{\kappa,N,\varepsilon,f,t}(x,F_1(x))^{\frac{c}{c+1}}
\left(\frac{J_1(x)}{\rho_0(x)}\right)^{\frac{c}{c+1}}
 \right)\mu_0(\d x)
\\
&= (1-t)\int_M\varphi_U\left(\left(
\frac{{\beta}_{\kappa,N,\varepsilon,f,1-t}(F_1(x),x)}{\rho_0(x)}
\right)^{\frac{c}{c+1}} \right)\mu_0(\d x)\\
&\hspace{2cm}
+t\int_M\varphi_U\left(\left(
\frac{{\beta}_{\kappa,N,\varepsilon,f,t}(x,F_1(x))}{\rho_1(F_1(x))}
\right)^{\frac{c}{c+1}} \right)\mu_0(\d x).
\end{align*}
Here we used the convexity and non-increasing property of $\varphi_U$ in the first inequality (cf. \cite[Remark 17.2]{Vi2}).
From $\pi=(F_0\times F_1)_{\sharp}\mu_0$, one can conclude the desired inequality.
\end{pf}

\subsection{Displacement convexity implies curvature bounds}\label{subsec:DisplacementconvCD}
The implication from \ref{item:Main2} to \ref{item:Main3} is trivial.
We now show that from \ref{item:Main3} to \ref{item:Main1},
and complete the proof of Theorem~\ref{thm:Main}.
For subsets $X,Y\subset M$ and $t\in[0,1]$, let $Z_t(X,Y)$ be the set of all points $\gamma(t)$, where $\gamma: [0, 1]\to  M$ is a minimal geodesic with $\gamma(0)\in X$, $\gamma(1)\in Y$.
We begin with the following Brunn-Minkowski inequality (cf. \cite[Lemma~4.3]{Sk:OneDim}):
\begin{lem}\label{lem:BM}
Let $X,Y\subset M$ be two bounded Borel subsets with $\m(X),\m(Y)\in]\,0, +\infty\,[$.
For $\kappa\in\R$, if $(M,d,\m)$ satisfies ${\rm TwCD_{rel}}(\kappa,N,\varepsilon)$,
then for every $t\in ]\,0,1\,[$,
\begin{align*}
\m(Z_t(X,Y))^{\frac{c}{c+1}}&\geq (1-t)\left(\inf_{(x,y)\in X\times Y}\beta_{\kappa,N,\varepsilon,f,1-t}(y,x)^{\frac{c}{c+1}} \right)\m(X)^{\frac{c}{c+1}}\nonumber\\
&\hspace{3cm}+t\left(\inf_{(x,y)\in X\times Y}\beta_{\kappa,N,\varepsilon,f,t}(x,y)^{\frac{c}{c+1}} \right)\m(Y)^{\frac{c}{c+1}}.
\end{align*}
\end{lem}
\begin{pf}
The proof is similar to that in \cite[Lemma~4.3]{Sk:OneDim}. 
We omit it.  
\end{pf}

Having Lemma~\ref{lem:BM} at hand, let us prove the following (cf. \cite[Proposition~4.4]{Sk:OneDim}):

\begin{prop}\label{prop:BM}
We suppose $\eps \neq 0$. 
For $\kappa\in\R$, if $(M,d,\m)$ satisfies ${\rm TwCD_{rel}}(\kappa,N,\varepsilon)$, then $\Ric_f^N\geq c^{-1}\kappa\,e^{-\frac{4(1-\varepsilon)f}{n-1}}g$. 
\end{prop}
\begin{pf}
We will follow the method of the proof of \cite[Theorem~1.2]{Oint}, \cite[Theorem~4.10]{Oh<0}.
Fix $x\in M$ and $v\in U_xM$,
and set
\begin{equation*}
 \theta_{\eps}:=-\frac{1}{n-1}\frac{1}{\eps}\left( \eps-\eps_0  \right)g(\nabla f,v).
\end{equation*}
Here we used the assumption $\varepsilon \neq 0$.
For a sufficiently small $t_0>0$,
let $\gamma:]-t_0,t_0[\to M$ be the geodesic with $\gamma(0)=x$ and $\dot{\gamma}(0)=v$.
Take $\delta\in ]0,t_0[$ and $\eta\in]0,\delta[$.
We denote by $B_{r}(o)$ the open geodesic ball of radius $r>0$ centered at $o\in M$,
and put $X:=B_{\eta(1+\theta_{\varepsilon} \delta)}(\gamma(-\delta))$ and $Y:=B_{\eta(1-\theta_{\varepsilon}\delta)}(\gamma(\delta))$.
Lemma~\ref{lem:BM} tells us that
\begin{align*}
\m\left(Z_{\frac12}(X,Y) \right)^{\frac{c}{c+1}}&\geq\frac12\left(\inf_{(x,y)\in X\times Y}\beta_{{\kappa,N,\varepsilon,f,\frac12}}(y,x)^{\frac{c}{c+1}} \right)\m(X)^{\frac{c}{c+1}}\\
&\hspace{1cm}+\frac12\left(\inf_{(x,y)\in X\times Y}\beta_{\kappa,N,\varepsilon,f,\frac12}(x,y)^{\frac{c}{c+1}} \right)\m(Y)^{\frac{c}{c+1}}.
\end{align*}
Letting $\eta\to0$ in the above inequality,
we have 
\begin{align}
\varliminf_{\eta\to0}\left(\frac{\m\left(Z_{\frac12}(X,Y) \right)}{\omega_n\eta^n} \right)^{\frac{c}{c+1}}&\geq 
\frac12\left(e^{-f(\gamma(-\delta))}(1+\theta_{\varepsilon} \delta)^n\beta_{\kappa,N,\varepsilon,f,\frac12}(\gamma(\delta),\gamma(-\delta)) \right)^{\frac{c}{c+1}}
\nonumber
\\&\hspace{1cm}+\frac12
\left(e^{-f(\gamma(\delta))}(1-\theta_{\varepsilon}\delta)^n\beta_{\kappa,N,\varepsilon,f,\frac12}(\gamma(-\delta),\gamma(\delta)) \right)^{\frac{c}{c+1}},\label{eq:BM1}
\end{align}
where $\omega_n$ is the volume of the unit ball in $\R^n$. 

Since
\begin{align}\label{eq:app1}
d_{N,\varepsilon,f,\frac12}(\gamma(\delta),\gamma(-\delta))&=\int_0^{\delta}e^{-\frac{2(1-\varepsilon)f(\gamma(\xi))}{n-1}}\d \xi,\\ \notag
d_{N,\varepsilon,f,\frac12}(\gamma(-\delta),\gamma(\delta))&=\int_{-\delta}^0e^{-\frac{2(1-\varepsilon)f(\gamma(\xi))}{n-1}}\d \xi,\\ \notag
d_{N,\varepsilon,f}(\gamma(\delta),\gamma(-\delta))&=\int_{-\delta}^{\delta}e^{-\frac{2(1-\varepsilon)f(\gamma(\xi))}{n-1}}\d \xi,
\end{align}
the Taylor series with respect to $\delta$ at 0 are
\begin{align}\label{eq:app2}
\beta_{\kappa,N,\varepsilon,f,\frac12}(\gamma(\delta),\gamma(-\delta))&=1-c^{-1}\frac{1-\varepsilon}{n-1}g(\nabla f,v)\delta\\ \notag
&+\left(c^{-1}\kappa e^{-\frac{4(1-\varepsilon)f(x)}{n-1}} +(1-c)\left(c^{-1}\frac{1-\varepsilon}{n-1} \right)^2 g(\nabla f,v)^2\right)\frac{\delta^2}{2}+O(\delta^3),\\ \notag
\beta_{\kappa,N,\varepsilon,f,\frac12}(\gamma(-\delta),\gamma(\delta))&=1+c^{-1}\frac{1-\varepsilon}{n-1}g(\nabla f,v)\delta\\ \notag
&+\left(    c^{-1}\kappa e^{-\frac{4(1-\varepsilon)f(x)}{n-1}} +(1-c)\left(c^{-1}\frac{1-\varepsilon}{n-1} \right)^2 g(\nabla f,v)^2   \right)\frac{\delta^2}{2}+O(\delta^3),
\end{align}
and
\begin{align}\label{eq:app3}
e^{-f(\gamma(-\delta))+f(x)}&=1+g(\nabla f,v) \delta+\left(g(\nabla f,v)^2-\nabla^2 f(v,v) \right)\frac{\delta^2}{2}+O(\delta^3),\\ \notag
e^{-f(\gamma(\delta))+f(x)}&=1-g(\nabla f,v)\delta+\left(g(\nabla f,v)^2-\nabla^2 f(v,v) \right)\frac{\delta^2}{2}+O(\delta^3),
\end{align}
and 
\begin{align}\label{eq:app4}
(1+\theta_{\varepsilon}\delta)^n&=1+n\theta_{\varepsilon}\delta+\frac{n(n-1)}{2}\theta^2_{\varepsilon}\delta^2+O(\delta^3),\\ \notag
(1-\theta_{\varepsilon}\delta)^n&=1-n\theta_{\varepsilon}\delta+\frac{n(n-1)}{2}\theta^2_{\varepsilon}\delta^2+O(\delta^3).
\end{align}
Substituting these series into \eqref{eq:BM1}, we have 
\begin{align}\label{eq:app6}
&\varliminf_{\eta\to0}\frac{\m\left(Z_{\frac12}(X,Y) \right)}{\omega_n\eta^n}\\
&\hspace{0.3cm}\geq
e^{-f(x)}\left\{1+\left(c^{-1}\kappa e^{-\frac{4(1-\varepsilon)f(x)}{n-1}}-\nabla^2 f(v,v)-
\frac{F(\theta_{\varepsilon})}{c+1} 
\right)\frac{\delta^2}{2}\right\}+O(\delta^3),\nonumber
\end{align}
where for $\alpha:=(1-\varepsilon)(n-1)^{-1}$ we set
\begin{equation}\label{eq:app5}
F(\theta):=n(n-(n-1)(c+1))\theta^2+2n(\alpha-c)g(\nabla f,v) \theta+(\alpha^2+2\alpha-c)g(\nabla f,v)^2.
\end{equation}
Now,
we can calculate
\begin{equation*}
F(\theta)+\frac{c+1}{N-n}g(\nabla f,v)^2=\frac{n}{\eps_0}\left\{\eps \theta+\frac{1}{n-1}(\eps-\eps_0)g(\nabla f,v)    \right\}^2,
\end{equation*}
and hence
\begin{align}\label{eq:Taylor}
&\varliminf_{\eta\to0}\frac{\m\left(Z_{\frac12}(X,Y) \right)}{\omega_n\eta^n}\\
&\hspace{0.3cm}\geq
e^{-f(x)}\left\{1+\left(c^{-1}\kappa e^{-\frac{4(1-\varepsilon)f(x)}{n-1}}-\nabla^2 f(v,v)+
\frac{g(\nabla f,v)^2}{N-n} 
\right)\frac{\delta^2}{2}\right\}+O(\delta^3)\nonumber
\end{align}
by the definition of $\theta_{\eps}$.
The detailed calculation can be seen in Appendix.

On the other hand, 
\begin{align}
\varlimsup_{\eta\to0}\frac{\m\left(Z_{\frac12}(X,Y) \right)}{\omega_n\eta^n}\leq e^{-f(x)}\left(1+\Ric_g(v)\frac{\delta^2}{2}\right)+O(\delta^3).\label{eq:RiciUpper}
\end{align} 
By comparing \eqref{eq:Taylor} and \eqref{eq:RiciUpper},
\begin{align*}
\Ric_g(v)&\geq c^{-1}\,\kappa\,e^{-\frac{4(1-\varepsilon)f(x)}{n-1}}-\nabla^2 f(v,v)+\frac{g(\nabla f,v)^2}{N-n},
\end{align*}
which means $\Ric_f^N(v)\geq c^{-1}\,\kappa\,e^{-\frac{4(1-\varepsilon)f(x)}{n-1}}$.
This completes the proof.
\end{pf}
We are now in a position to conclude Theorem~\ref{thm:Main}.

\begin{apf}{Theorem~\ref{thm:Main}}
By Propositions \ref{prop:CDimpliesConvexity} and \ref{prop:BM},
we complete the proof.
\end{apf}

\subsection{Interpolation inequalities}\label{subsec:Inteplation} 
Under the curvature condition \eqref{eq:RicciLowerBddWY},
the second named author \cite{Sk:OneDim} has derived some interpolation inequalities from the proof of the characterization result (see \cite[Subsection 4.3]{Sk:OneDim}).
By the same argument,
we can obtain such interpolation inequalities in our setting,
and we collect them here.
We just present their forms,
and the proof is left to the readers.

We start with the $p$-mean inequality (cf. \cite[Corollary 4.5]{Sk:OneDim}). Let $t\in]\,0, 1\,[$ and $a,b\in[\,0,+\infty\,[$. For $p\in\R\setminus\{0\}$, the {\it $p$-mean} is defined as follows:
\begin{align*}
\mathcal{M}_t^p(a,b):=((1-t)a^p+tb^p)^{\frac{1}{p}}
\end{align*}
if $ab\ne0$, and $\mathcal{M}_t^p(a,b):=0$ if $ab=0$. As the limits, it is defined as 
\begin{align*}
\mathcal{M}_t^0(a,b):=a^{1-t}b^t,\quad \mathcal{M}_t^{\infty}(a,b):=\max\{a,b\},\quad 
\mathcal{M}_t^{-\infty}(a,b):=\min\{a,b\}.
\end{align*}

\begin{cor}\label{cor:PrekopaLeindler1}
For $i=0,1$, let $\psi_i:M\to\R$ be non-negative, integrable functions. Let $X,Y\subset M$ be bounded Borel subsets with ${\rm supp}\,[\psi_0]\subset X$, ${\rm supp}\,[\psi_1]\subset Y$. Let $\psi:M\to\R$ be a non-negative function. For $t\in]\,0,1\,[$ and $p\geq -c(c+1)^{-1}$, we assume that for all $(x, y)\in X\times Y$ and $z\in Z_t(\{x\}, \{y\})$, we have 
\begin{align*}
\psi(z)\geq\mathcal{M}_t^p\left(\frac{\psi_0(x)}{\beta_{\kappa,N,\varepsilon,f,1-t}(y,x)},\frac{\psi_1(y)}{\beta_{\kappa,N,\varepsilon,f,t}(x,y)} \right).
\end{align*}
For $\kappa\in\R$, if \,$\Ric_f^N\geq c^{-1}\,\kappa\,e^{-\frac{4(1-\varepsilon)f}{n-1}}g$, then we have 
\begin{align*}
\int_M\psi\,\d\m\geq\mathcal{M}_t^{\frac{cp}{(1+c)p+c}}
\left(\int_M\psi_0\,\d\m,\int_M\psi_1\,\d\m \right).
\end{align*}
Here we set $c\,p((1+c)p+c)^{-1}:=-\infty$ for $p=-c(c+1)^{-1}$.
\end{cor} 

We next show the Pr\'ekopa-Leindler inequality, which is the case of $p=0$ in Corollary \ref{cor:PrekopaLeindler1} (cf. \cite[Corollary 4.6]{Sk:OneDim}): 

\begin{cor}\label{cor:PrekopaLeindler2}
For $i=0,1$, let $\psi_i, X, Y, \psi$ be as in 
Corollary~\ref{cor:PrekopaLeindler1}. For $t\in]\,0,1\,[$, 
we assume that for all $(x,y)\in X \times Y$ and $z\in  Z_t(\{x\},\{y\})$,
\begin{align*}
\psi(z)\geq \left(\frac{\psi_0(x)}{\beta_{\kappa,N,\varepsilon,f,1-t}(y,x)} \right)^{1-t}\left(\frac{\psi_1(y)}{\beta_{\kappa,N,\varepsilon,f,t}(x,y)} \right)^t.
\end{align*}
For $\kappa\in\R$, if \,$\Ric_f^N\geq c^{-1}\,\kappa\,e^{-\frac{4(1-\varepsilon)f}{n-1}}g$, then we have 
\begin{align*}
\int_M\psi\,\d\m\geq
\left(\int_M\psi_0\d\m\right)^{1-t}\left(\int_M\psi_1\,\d\m \right)^t.
\end{align*}
\end{cor}

We further possess the Borell-Brascamp-Lieb inequality, which is the case of $p=-c(c+1)^{-1}$ in Corollary \ref{cor:PrekopaLeindler1} (cf. \cite[Corollary 4.7]{Sk:OneDim}): 

\begin{cor}\label{cor:PrekopaLeindler3}
For $i=0,1$, let $\psi_i, X, Y, \psi$ be as in Corollary ~\ref{cor:PrekopaLeindler1}. We suppose $\int_M\psi_0\,\d\m=
\int_M\psi_1\,\d\m=1$. For $t\in]\,0,1\,[$, we assume that for all $(x,y)\in X\times Y$ and $z\in Z_t(\{x\},\{y\})$,
\begin{align*}
\psi(z)^{-\frac{c}{c+1}}\leq (1-t)\left(\frac{\psi_0(x)}{\beta_{\kappa,N,\varepsilon,f,1-t}(y,x)} \right)^{-\frac{c}{c+1}}+t\left(\frac{\psi_1(y)}{\beta_{\kappa,N,\varepsilon,f,t}(x,y)} \right)^{-\frac{c}{c+1}}.
\end{align*}
For $\kappa\in\R$, if \,$\Ric_f^N\geq c^{-1}\,\kappa\,e^{-\frac{4(1-\varepsilon)f}{n-1}}g$, then we have 
$\int_M\psi\,\d\m\geq1$. 
\end{cor}

%

\section{Functional Inequalities}\label{sec:FuncIneq}
In this last section,
we discuss functional inequalities under the curvature condition \eqref{eq:WeightedLowerBddd}.
For $\kappa\in\R$, let 
$\mathfrak{c}_{\kappa}:=\s_{\kappa}'$. 
Following \cite[Section 5]{Sk:OneDim},
for $x,y\in M$
we define  
\begin{align*}
{\sf b}_{\kappa,N,\varepsilon,f}(x,y)&:=\left(\frac{e^{-\frac{2(1-\varepsilon)f(x)}{n-1}}d(x,y)}{\s_{\kappa}(d_{N,\varepsilon,f}(x,y))} \right)^{c^{-1}},\\ 
\mathfrak{b}_{\kappa,N,\varepsilon,f}(x,y)&:=
\frac{1}{c+1}\left(\frac{e^{-\frac{2(1-\varepsilon)f(x)}{n-1}}d(x,y)\mathfrak{c}_{\kappa}(d_{N,\varepsilon,f}(x,y))}{\s_{\kappa}(d_{N,\varepsilon,f}(x,y))}-1 \right)
\end{align*}
if $d_{N,\varepsilon,f}(x,y)\in]0,C_{\kappa}[$;
${\sf b}_{\kappa,N,\varepsilon,f}(x,y):=1$ and $\mathfrak{b}_{\kappa,N,\varepsilon,f}(x,y):=0$ if $x=y$;
otherwise, ${\sf b}_{\kappa,N,\varepsilon,f}(x,y):=+\infty$ and 
$\mathfrak{b}_{\kappa,N,\varepsilon,f}(x,y):=+\infty$ (cf.~Remark~\ref{rem:TwistedCoeffic}). 

One can verify the following fact (cf. \cite[Lemma 5.2]{Sk:OneDim}). The proof is left to the readers.

\begin{lem}\label{lem:fucIneq}
Let $\kappa\in\R$.
Let $x,y\in M$ satisfy $d_{N,\varepsilon,f}(x,y)\in [0,C_{\kappa}[$. If $y\notin {\rm Cut}(x)$,
then as $t\to0$, we have
\begin{align*}
\beta_{\kappa,N,\varepsilon,f,t}(x,y)\to{\sf b}_{\kappa,N,\varepsilon,f}(x,y),\quad \frac{1-\beta_{\kappa,N,\varepsilon,f,1-t}(y,x)^{\frac{c}{c+1}}}{t}\to\mathfrak{b}_{\kappa,N,\varepsilon,f}(x,y).
\end{align*}
\end{lem}

For a non-negative Lipschitz function $\rho$ on $M$ with 
$\int_M\rho\,\d\m=1$, set $\mu:=\rho\,\m$. 
The \emph{generalized Fisher information} $I_{\m}(\mu)$ of $\mu$ is defined as 
\begin{align*}
I_{\m}(\mu):=\int_M\frac{\|\nabla\rho^{\frac{1}{c+1}}\|^2}{\rho}\d\m.
\end{align*}
In general, $I_{\m}(\mu)\in[\,0,+\infty\,]$.
%
%
We present the following (cf. \cite[Proposition 5.4]{Sk:OneDim}):
\begin{prop}\label{prop:Fisher2}
Suppose that $\m(M)<+\infty$ and $\m$ has finite second moment.
For $i=0,1$, let $\rho_i:M\to\R$ be non-negative Lipschitz functions with $\int_M\rho_i\,\d\m=1$.
We assume that
$\mu:=\rho_0\,\m$ and $\nu:=\rho_1\,\m$ belong to $\mathscr{P}_2^{\rm ac}(M)$.
For $\kappa\in\R$, if ${\Ric}_f^N\geq c^{-1}\kappa \,e^{-\frac{4(1-\varepsilon)f}{n-1}}g$, then
\begin{align*}
H_{\m}(\mu)&\leq\sqrt{I_{\m}(\mu)}W_2(\mu,\nu)+
\frac{c+1}{c}\int_{M^2}\rho_0(x)^{-\frac{c}{c+1}}
\mathfrak{b}_{\kappa,N,\varepsilon,f}(x,y)\pi(\d x\d y)
\\
&\hspace{4cm}- \frac{c+1}{c}\int_{M^2}\rho_1(y)^{-\frac{c}{c+1}}
\left(
{\sf b}_{\kappa,N,\varepsilon,f}(x,y)^{\frac{c}{c+1}}-1
\right)\pi(\d x\d y)\nonumber\\
&\hspace{5cm}- \frac{c+1}{c}\int_{M^2}\left(\rho_1(y)^{-\frac{c}{c+1}}-1 \right)\pi(\d x \d y),\nonumber
\end{align*}
where $\pi$ is a unique optimal coupling of $(\mu,\nu)$ with respect to the square of distance. 
When $I_{\m}(\mu)=+\infty$ and $\mu=\nu$, we use the convention 
$I_{\m}(\mu)W_2(\mu,\mu)=0$.
\end{prop}
\begin{pf}
If $I_{\m}(\mu)=+\infty$ and $\mu\ne\nu$, the inequality trivially holds. 
We first assume $I_{\m}(\mu)<+\infty$. 
By virtue of Corollary~\ref{cor:Main}, $(M, d, \m)$ satisfies ${\rm TwCD_{rel}}(\kappa,N,\varepsilon)$,
and hence
\begin{align*}
H_{\m}(\mu_t)&\leq\frac{c+1}{c}-\frac{c+1}{c}(1-t)\int_{M^2}\rho_0(x)^{-\frac{c}{c+1}}\beta_{\kappa,N,\varepsilon,f,1-t}(y,x)^{\frac{c}{c+1}}\pi(\d x \d y)\\
&\hspace{4cm}-
\frac{c+1}{c} t\int_{M^2}\rho_1(y)^{-\frac{c}{c+1}}\beta_{\kappa,N,\varepsilon,f,t}(x,y)^{\frac{c}{c+1}}\pi(\d x \d y),
\end{align*}
here $(\mu_t)_{t\in[0,1]}$ is a unique minimal geodesic in $(\mathscr{P}_2(M),W_2)$ from $\mu$ to $\nu$. Therefore,
\begin{align*}
\frac{H_{\m}(\mu_t)-H_{\m}(\mu)}{t}&\leq 
\frac{c+1}{c}\int_{M^2}\rho_0(x)^{-\frac{c}{c+1}}
\frac{1-\beta_{\kappa,N,\varepsilon,f,1-t}(y,x)^{\frac{c}{c+1}}}{t}\pi(\d x\d y)\\
&\hspace{0.5cm}+
\frac{c+1}{c}\int_{M^2}\rho_0(x)^{-\frac{c}{c+1}}
\left(\beta_{\kappa,N,\varepsilon,f,1-t}(y,x)^{\frac{c}{c+1}}-1 \right)\pi(\d x\d y)\\
&\hspace{0.5cm}-
\frac{c+1}{c}\int_{M^2}\rho_1(y)^{-\frac{c}{c+1}}
\left(\beta_{\kappa,N,\varepsilon,f,t}(x,y)^{\frac{c}{c+1}}-1 \right)\pi(\d x\d y)\\
&\hspace{0.5cm}-
\frac{c+1}{c}\int_{M^2}\left(\rho_1(y)^{-\frac{c}{c+1}}-1\right)\pi(\d x\d y)-H_{\m}(\mu).
\end{align*}
Let $F_1$ be the map defined as \eqref{eq;OptimalMap}.
We can deduce $d_{N,\varepsilon,f}(x,F_1(x))\in [0,C_{\kappa}[$ for $\mu$-almost every $x\in M$ from Theorem~\ref{thm:CorderoMcCannSche} and Proposition~\ref{prop:DirectionalMyers}.
Lemma~\ref{lem:fucIneq} and $\pi=(F_0\times F_1)_{\sharp}\mu$ yields
\begin{align*}
\varlimsup_{t\to0}\frac{H_{\m}(\mu_t)-H_{\m}(\mu)}{t}&\leq 
\frac{c+1}{c}\int_{M^2}\rho_0(x)^{-\frac{c}{c+1}}
\mathfrak{b}_{\kappa,N,\varepsilon,f}(x,y)
\pi(\d x\d y)\\
&\hspace{0.5cm}-
\frac{c+1}{c}\int_{M^2}\rho_1(y)^{-\frac{c}{c+1}}
\left({\sf b}_{\kappa,N,\varepsilon,f}(x,y)^{\frac{c}{c+1}}-1 \right)\pi(\d x\d y)\\
&\hspace{0.5cm}-
\frac{c+1}{c}\int_{M^2}\left(\rho_1(y)^{-\frac{c}{c+1}}-1\right)\pi(\d x\d y)-H_{\m}(\mu).
\end{align*}
Now we assume that $\rho_0$ is bounded below away from $0$. Then one can apply 
\cite[Theorem~20.1]{Vi2} so that
\begin{align}
\varliminf_{t\to0}\frac{H_{\m}(\mu_t)-H_{\m}(\mu)}{t}\geq 
-\sqrt{I_{\m}(\mu)}W_2(\mu,\nu)\label{eq:FisherIneq}
\end{align}
under $I_{\m}(\mu)<+\infty$ (cf. \cite[Remark 20.2]{Vi2}). 
Note that the condition $\m(M)<+\infty$ assures the integrability conditions 
$H(\rho_0),\rho_0H'(\rho_0)\in L^1(M;\m)$ in \cite[Theorem~20.1]{Vi2}, because 
\begin{align*}
&\quad \,\,\int_M\rho_0(x)^{\frac{1}{c+1}}\m(\d x)\\
&\leq \left(\int_M(1+d(x,x_0)^2)\rho_0(x)\m(\d x) \right)^{\frac{1}{c+1}}\left(\int_M\frac{\m(\d x)}{(1+d(x,x_0)^2)^{1/c}} \right)^{\frac{c}{c+1}}<+\infty.
\end{align*} 
Then we have the conclusion. Next we prove the assertion without assuming the existence of positive lower bound for $\rho_0$. Now we set 
\begin{align*}
\rho_0^i:=\frac{i\rho_0+1}{i+\m(M)}\quad\text{ and }\quad \mu^i:=\rho_0^i\m\in\mathscr{P}_2^{\rm ac}(M). 
\end{align*}
Then we have the conclusion by replacing $\mu$ (resp.~$\rho_0$) with $\mu^i$ 
(resp.~$\rho_0^i$). Since $\m$ has finite second moment, we see $W_2(\mu^i,\mu)\to0$ as $i\to +\infty$. So the conclusion can be obtained under $I_{\m}(\mu)<+\infty$, 
because of the lower semi continuity of $\nu\mapsto H_{\m}(\nu)$ in $(\mathscr{P}_2(M),W_2)$. If $\mu=\nu$, then $\mu_t=\mu=\nu$, hence \eqref{eq:FisherIneq} also trivially holds  even if $I_{\m}(\mu)=+\infty$.
\end{pf}


We will show three functional inequalities under 
the curvature condition \eqref{eq:WeightedLowerBddd}.
In what follows,
we always assume $\m \in \mathscr{P}_2^{\rm ac}(M)$.
To state our results,
we introduce the following condition,
which seems to be quite strong:
We say that $\mu\in\mathscr{P}_2^{\rm ac}(M)$ is {\it $\m$-constant} if $d_{N,\eps,f}(x,F_1(x))=e^{-\frac{2(1-\varepsilon)f(x)}{n-1}}d(x,F_1(x))$ on $M$,
where $F_1$ is the map defined as \eqref{eq;OptimalMap} for $\nu=\m$.
We obtain the following (cf. \cite[Theorems 20.10, 21.7]{Vi2}):

\begin{cor}\label{cor:HWI-LSI}
We assume $\m \in \mathscr{P}_2^{\rm ac}
(M)$.
Let $\rho:M\to\R$ be a 
non-negative Lipschitz function with $\int_M\rho\,\d\m=1$.
Assume that $\mu:=\rho\,\m$ belongs to $\mathscr{P}_2^{\rm ac}(M)$. We further assume that $(1-\eps)f\leq (n-1)\delta$ for $\delta \in \mathbb{R}$,
and $\mu$ is $\m$-constant. For $\kappa>0$, if \,$\Ric_f^N\geq c^{-1}\,\kappa\,e^{-\frac{4(1-\varepsilon)}{n-1}f}g$, then we have  
\begin{enumerate}
\item The HWI inequality 
\begin{align*}
H_{\m}(\mu)\leq\sqrt{I_{\m}(\mu)}W_2(\mu,\m)-\frac{\kappa e^{-4\delta}}{6\,c}\left(1+2(\sup\rho)^{-\frac{c}{c+1}} \right)W_2(\mu,\m)^2;
\end{align*}
\item the Logarithmic Sobolev inequality
\begin{align*}
H_{\m}(\mu)\leq \frac{3c\left(1+2(\sup\rho)^{-\frac{c}{c+1}}  \right)^{-1}}{2\kappa e^{-4\delta}}I_{\m}(\mu).
\end{align*}
\end{enumerate}
\end{cor}
\begin{pf}
Note first that $M$ is compact under $\kappa>0$ and $(1-\eps)f\leq (n-1)\delta$ for $\delta \in \mathbb{R}$ (see \cite[Proposition~3.2]{KwSk}). 
Hence $\rho$ is bounded.
We begin with the HWI inequality. Since $\pi=(F_0\times F_1)_{\sharp}\mu$, and $\mu$ is $\m$-constant,
\begin{align*}
{\sf b}_{\kappa,N,\eps,f}(x,y)&=\left(\frac{d_{N,\eps,f}(x,y)}{\s_{\kappa}(d_{N,\eps,f}(x,y))} \right)^{c^{-1}},\\
\mathfrak{b}_{\kappa,N,\eps,f}(x,y)&=\frac{1}{c+1}\left(\frac{d_{N,\eps,f}(x,y)\mathfrak{c}_{\kappa}(d_{N,\eps,f}(x,y))}{\s_{\kappa}(d_{N,\eps,f}(x,y))}-1 \right)
\end{align*}
on the support of $\pi$. By elementary estimates and $(1-\eps)f\leq (n-1)\,\delta$,
we possess the following (cf. \cite[(20.32), (20.34)]{Vi2}):
\begin{align*}
-\left({\sf b}_{\kappa,N,\eps,f}(x,y)^{\frac{c}{c+1}}-1 \right)&\leq -\frac{\kappa}{6(c+1)}d_{N,\eps,f}(x,y)^2\leq -\frac{\kappa e^{-4\delta}}{6(c+1)}d(x,y)^2,\\
\mathfrak{b}_{\kappa,N,\eps,f}(x,y)&\leq- \frac{\kappa}{3(c+1)}d_{N,\eps,f}(x,y)^2\leq  -\frac{\kappa e^{-4\delta}}{3(c+1)}d(x,y)^2.
\end{align*}
Applying Proposition~\ref{prop:Fisher2} to $\rho_0 =\rho$ and $\rho_1 = 1$, we see
\begin{align*}
H_{\m}(\mu)&\leq\sqrt{I_{\m}(\mu)}W_2(\mu,\m)-
\frac{\kappa e^{-4\delta}}{3c}
\int_{M^2}\rho(x)^{-\frac{c}{c+1}}d(x,y)^2
\pi(\d x\d y)\\
&\hspace{2cm}-\frac{\kappa e^{-4\delta}}{6c}
\int_{M^2}d(x,y)^2\pi(\d x \d y),
\end{align*}
and hence
\begin{align*}
H_{\m}(\mu)\leq\sqrt{I_{\m}(\mu)}W_2(\mu,\m)-
\frac{\kappa e^{-4\delta}}{6c}
\left(1+2(\sup\rho)^{-\frac{c}{c+1}} \right)
\int_{M^2}d(x,y)^2\pi(\d x\d y).
\end{align*}
By the optimality of $\pi$, the right hand side of the above inequality is equal to that of the desired one. 
We next show the Logarithmic Sobolev inequality.
Using an elementary inequality, we have
\begin{align*}
\sqrt{I_{\m}(\mu)}W_2(\mu,\m)&\leq \frac{3c\,\left(1+2(\sup\rho)^{-\frac{c}{c+1}}  \right)^{-1}}{2\kappa e^{-4\delta}}I_{\m}(\mu)+\frac{\kappa\,e^{-4\delta}}{6c}\left( 1+2(\sup\rho)^{-\frac{c}{c+1}} \right)W_2(\mu,\m)^2.
\end{align*}
From the HWI inequality, one can derive the desired one. This completes the proof.
\end{pf}

Finally, we conclude the following finite dimensional transport energy inequality (cf. \cite[Theorem 22.37, Corollary 22.39]{Vi2}):

\begin{cor}\label{cor:TEI}
We assume $\m \in \mathscr{P}_2^{\rm ac}(M)$.
Let $\rho : M\to \R$ be a non-negative Lipschitz function with $\int_M\rho\,\d\m=1$.
Assume that $\mu:=\rho\,\m$ belongs to $\mathscr{P}_2^{\rm ac}(M)$,
and also assume that $\mu$ is $\m$-constant. 
For $\kappa>0$, if \,$\Ric_f^N\geq c^{-1}\kappa e^{-\frac{4(1-\varepsilon)f}{n-1}}g$, then 
\begin{align*}
H_{\m}(\mu)&\geq \frac12\cdot\frac{c+1}{c}
+\frac12\int_M\rho^{\frac{1}{c+1}}\log\rho\,\d\m\\
&\hspace{2cm}-\frac12\cdot\frac{c+1}{c}\int_{M^2}\left(\mathfrak{b}_{\kappa,N,\eps,f}(x,y)+\exp\left(1-{\sf b}_{\kappa,N,\eps,f}(x,y)^{\frac{c}{c+1}} \right) \right)\pi(\d x\d y),
\end{align*}
where $\pi$ is the unique optimal coupling of $(\mu,\m)$. 
\end{cor}
\begin{pf}
Under $\m\in \mathscr{P}_2^{\rm ac}(M)$, we see the well-definedness of $\int_M\rho^{\frac{1}{c+1}}\log\rho\,\d\m\in]-\infty,+\infty\,]$, because $x^{a}\log x$ is bounded below for any $a \in]0,1]$. 
We start with 
\begin{align}
2H_{\m}(\mu)=2\cdot\frac{c+1}{c}-2\cdot\frac{c+1}{c}\int_{M^2}\rho(x)^{-\frac{c}{c+1}}\pi(\d x \d y).\label{eq:Hm1}
\end{align}
Let us recall the following Young inequality:
\begin{align*}
ab\leq a\log a-2a+e^{b+1}.
\end{align*}
We set $\mathcal{B}(x,y):={\sf b}_{\kappa,N,\eps,f}(x,y)^{\frac{c}{c+1}}$. From the Young inequality, we derive 
\begin{align*}
\rho(x)^{-\frac{c}{c+1}}\log\rho(x)^{\frac{c}{c+1}}&=\left(\rho(x)^{-\frac{c}{c+1}}e^{-\mathcal{B}(x,y)} \right)\left(e^{\mathcal{B}(x,y)}\log\rho(x)^{\frac{c}{c+1}} \right)\\
&\leq \left(\rho(x)^{-\frac{c}{c+1}}e^{-\mathcal{B}(x,y)}  \right)
\left(e^{\mathcal{B}(x,y)}\mathcal{B}(x,y)-2e^{\mathcal{B}(x,y)}+e\rho(x)^{\frac{c}{c+1}} \right)\\
&=\rho(x)^{-\frac{c}{c+1}}\mathcal{B}(x,y)-2\rho(x)^{-\frac{c}{c+1}}+e^{1-\mathcal{B}(x,y)}
\end{align*}
on the support of $\pi$, and hence
\begin{align}
-2\rho(x)^{-\frac{c}{c+1}}\geq \rho(x)^{-\frac{c}{c+1}}\log\rho(x)^{\frac{c}{c+1}}
-\rho(x)^{-\frac{c}{c+1}}\mathcal{B}(x,y)-e^{1-\mathcal{B}(x,y)}.\label{eq:Hm2}
\end{align}
By \eqref{eq:Hm1} and \eqref{eq:Hm2}, we obtain
\begin{align}\label{eq:Conclusion1}
2H_{\m}(\mu)&\geq 2\cdot\frac{c+1}{c}+\int_M\rho^{\frac{1}{c+1}}\log\rho\,\d\m\\ \notag
&\hspace{1cm}- \frac{c+1}{c}\int_{M^2}\left(\rho(x)^{-\frac{c}{c+1}}
{\sf b}_{\kappa,N,\eps,f}(x,y)^{\frac{c}{c+1}}+\exp\left( 1-{\sf b}_{\kappa,N,\eps,f}(x,y)^{\frac{c}{c+1}}\right) \right)\pi(\d x\d y).
\end{align}
We apply Proposition~\ref{prop:Fisher2} to $\rho_0=1$ and $\rho_1=\rho$. 
From $H_{\m}(\m)=0$ and $I_{\m}(\m)=0$,
\begin{align*}
0\leq \int_{M^2}\left(\mathfrak{b}_{\kappa,N,\eps,f}(x,y)-\rho(y)^{-\frac{c}{c+1}}{\sf b}_{\kappa,N,\eps,f}(x,y)^{\frac{c}{c+1}} \right)\wh{\pi}(\d x\d y)+1,
\end{align*} 
 where $\wh{\pi}$ is a unique optimal coupling of $(\m,\mu)$. 
Since $\mu$ is $\m$-constant, $\mathfrak{b}_{\kappa,N,\eps,f}$  and ${\sf b}_{\kappa,N,\eps,f}$ are symmetric on the support of $\pi$. It follows that
\begin{align*}
0\leq \int_{M^2}\left(\mathfrak{b}_{\kappa,N,\eps,f}(x,y)-\rho(x)^{-\frac{c}{c+1}}{\sf b}_{\kappa,N,\eps,f}(x,y)^{\frac{c}{c+1}} \right){\pi}(\d x\d y)+1,
\end{align*}
which is equivalent to 
\begin{align}
-\int_{M^2}\rho(x)^{ -\frac{c}{c+1}}{\sf b}_{\kappa,N,\eps,f}(x,y)^{\frac{c}{c+1}}\pi(\d x\d y)\geq -1-\int_{M^2}
\mathfrak{b}_{\kappa,N,\eps,f}(x,y)\pi(\d x \d y).\label{eq:Conclusion2}
\end{align}
Combining \eqref{eq:Conclusion1} and \eqref{eq:Conclusion2} leads to the desired inequality.
\end{pf}
On Corollaries \ref{cor:HWI-LSI} and \ref{cor:TEI},
the authors do not know whether the assumption that $(1-\eps)f\leq (n-1)\delta$ and $\mu$ is $\m$-constant can be dropped.

Under the curvature condition \eqref{eq:RicciLowerBdd}, similar functional inequalities are known to be useful to analyze the gradient flow of entropy functionals (see e.g., \cite[Chapters 23, 24, 25]{Vi2}). There might be some applications of our inequalities to the analysis of such gradient flow under the curvature condition \eqref{eq:WeightedLowerBddd}.

Concerning the curvature condition \eqref{eq:RicciLowerBdd},
functional inequalities can be also derived from the so-called Bakry-\'Emery's $\Gamma$-calculus (\cite{BE1}).
Li-Xia \cite{LiXia} have formulated a Bochner type formula that is associated with $\Ric^1_f$.
One might be able to develop the $\Gamma$-calculus in our framework via their Bochner formula.

\section{Appendix}\label{sec:Appendix}

\subsection{Twisted coefficients}\label{subsec:twisted}
In this appendix,
we give a proof of the assertion stated in Remark \ref{rem:TwistedCoeffic}.
Namely, we show:
\begin{prop}
For the twisted coefficient $\beta_{\kappa,N,\varepsilon,f,t}(x,y)$ defined as \eqref{eq:twistedcoeffi},
it holds that
$\beta_{\kappa,N,\varepsilon,f,t}(x,y)\to 1$ as $d(x,y)\to 0$.
\end{prop}
\begin{pf}
It suffices to prove that
\begin{equation*}
\frac{d_{N,\varepsilon,f,t}(x,y)}{td_{N,\varepsilon,f}(x,y)}\to 1
\end{equation*}
as $d(x,y)\to 0$.
Fix $x\in M$, and a sufficiently small $r>0$.
Take $y\in B_r(x)$,
and a unique minimal geodesic $\gamma:[0,d(x,y)]\to M$ from $x$ to $y$.
We set $a:=2(n-1)^{-1}(1-\varepsilon)$.
Then it holds that
\begin{align*}
\left|\frac{d_{N,\varepsilon,f,t}(x,y)}{td_{N,\varepsilon,f}(x,y)}-1\right|&=
\left|  
\frac{\displaystyle \int_0^{td(x,y)}e^{-a f(\gamma(\xi))}\d\xi }{t \displaystyle\int_0^{d(x,y)}e^{-a f(\gamma(\xi))}\d\xi}
-1
\right|\\
&=
\frac{1}{\displaystyle \int_0^{d(x,y)}e^{-a f(\gamma(\xi))}\d\xi}
\left|
\int_0^{d(x,y) }e^{-a f(\gamma(t\xi))}\d\xi
-
\int_0^{d(x,y) }e^{-a f(\gamma(\xi))}\d\xi
\right|
\\
&\leq
\frac{1}{\displaystyle\int_0^{d(x,y)}e^{-a f(\gamma(\xi))}\d\xi}
\int_0^{d(x,y)}
e^{-a f(\gamma(\xi))}
\left| 1-e^{-a ( f(\gamma(t\xi))-f(\gamma(\xi))) }
\right|\d\xi.
\end{align*}
We now recall the following elementary estimate: For all $b\in \mathbb{R}$,
\begin{equation*}
\left|1-e^{-b} \right|\leq e^{\left|b \right|}-1.
\end{equation*}
It follows that
\begin{equation*}
\left| 1-e^{-a ( f(\gamma(t\xi))-f(\gamma(\xi))) }\right|\leq e^{|a| |f(\gamma(t\xi))-f(\gamma(\xi))| }-1.
\end{equation*}
Furthermore,
setting
\begin{equation*}
A:=(1-t)\sup_{B_r(x)}\|\nabla f\|,
\end{equation*}
we obtain
\begin{equation*}
|f(\gamma(t\xi))-f(\gamma(\xi))| \leq d(\gamma(t\xi), \gamma(\xi)) \sup_{B_r(x)}\|\nabla f\| \leq Ad(x,y). 
\end{equation*}
Therefore, we see
\begin{align*}
e^{|a| |f(\gamma(t\xi))-f(\gamma(\xi))| }-1&\leq e^{|a| Ad(x,y)}-1
= \sum_{k=1}^{\infty}\frac{ (|a|Ad(x,y))^k}{k!}\\
&=|a| Ad(x,y) \sum_{k=1}^{\infty}\frac{ (|a|Ad(x,y))^{k-1}}{k!}\\
&\leq |a|Ad(x,y) \sum_{k=1}^{\infty}\frac{ (|a|Ar)^{k-1}}{(k-1)!}=|a|Ae^{|a|Ar}d(x,y).
\end{align*}
Combining the above estimates,
we arrive at
\begin{equation*}
\left|\frac{d_{N,\varepsilon,f,t}(x,y)}{td_{N,\varepsilon,f}(x,y)}-1\right|\leq |a|A e^{|a|Ar}d(x,y).
\end{equation*}
This proves the desired claim.
\end{pf}

\subsection{Taylor series}\label{subsec:taylor}

This appendix is also devoted to a supplemental material for the proof of Proposition \ref{prop:BM} since the calculation is straightforward but quite complicated.
We use the same notation as in the proof.

First,
we give an outline of the proof of \eqref{eq:app2}.
In view of \eqref{eq:app1},
one can verify
\begin{align*}
d_{N,\varepsilon,f,\frac{1}{2}}(\gamma(\delta),\gamma(-\delta))&=e^{\frac{-2(1-\eps)f(x)}{n-1}}\delta\left( 1-\frac{(1-\eps)g(\nabla f,v)}{n-1}\delta+A\,\delta^2 +O(\delta^3)   \right),\\ 
d_{N,\varepsilon,f,\frac{1}{2}}(\gamma(-\delta),\gamma(\delta))&=e^{\frac{-2(1-\eps)f(x)}{n-1}}\delta\left( 1+\frac{(1-\eps)g(\nabla f,v)}{n-1}\delta+A \,\delta^2 +O(\delta^3)   \right),\\
d_{N,\varepsilon,f}(\gamma(\delta),\gamma(-\delta))&=2e^{\frac{-2(1-\eps)f(x)}{n-1}}\delta\left( 1+A \,\delta^2 +O(\delta^3)   \right),
\end{align*}
where
\begin{equation*}
A:=\frac{2(1-\eps)^2g(\nabla f,v)^2}{3(n-1)^2}-\frac{(1-\eps)\nabla^2 f(v,v)}{3(n-1)}.
\end{equation*}
Using $\mathfrak{s}_{\kappa}(s)=s-(\kappa/6)s^3+O(s^5)$,
we see
\begin{align*}
&\quad\,\, \mathfrak{s}_{\kappa}(d_{N,\varepsilon,f,\frac{1}{2}}(\gamma(\delta),\gamma(-\delta)))\\ \notag
&=e^{\frac{-2(1-\eps)f(x)}{n-1}}\delta\left\{ 1-\frac{(1-\eps)g(\nabla f,v)}{n-1}\delta+\left(A-\frac{\kappa}{6}e^{\frac{-4(1-\eps)f(x)}{n-1}} \right)\delta^2 +O(\delta^3)   \right\},\\
&\quad \,\,\mathfrak{s}_{\kappa}(d_{N,\varepsilon,f,\frac{1}{2}}(\gamma(-\delta),\gamma(\delta)))\\ \notag
&=e^{\frac{-2(1-\eps)f(x)}{n-1}}\delta\left\{ 1+\frac{(1-\eps)g(\nabla f,v)}{n-1}\delta+\left(A-\frac{\kappa}{6}e^{\frac{-4(1-\eps)f(x)}{n-1}} \right)\delta^2 +O(\delta^3)   \right\},\\
&\quad \,\,\mathfrak{s}_{\kappa}(d_{N,\varepsilon,f}(\gamma(\delta),\gamma(-\delta)))\\ \notag
&=2e^{\frac{-2(1-\eps)f(x)}{n-1}}\delta\left\{ 1+\left(A-\frac{2\kappa}{3}e^{\frac{-4(1-\eps)f(x)}{n-1}} \right)\delta^2 +O(\delta^3)   \right\}.
\end{align*}
By $(1+s)^{-1}=1-s+s^2-s^3+O(s^4)$,
we obtain
\begin{align*}
\frac{2\mathfrak{s}_{\kappa}(d_{N,\varepsilon,f,\frac{1}{2}}(\gamma(\delta),\gamma(-\delta)))}{\mathfrak{s}_{\kappa}(d_{N,\varepsilon,f}(\gamma(\delta),\gamma(-\delta)))}&=1-\frac{(1-\eps)g(\nabla f,v)}{n-1}\delta+\frac{\kappa}{2}e^{\frac{-4(1-\eps)f(x)}{n-1}}\delta^2+O(\delta^3),\\ \label{eq:7th}
\frac{2\mathfrak{s}_{\kappa}(d_{N,\varepsilon,f,\frac{1}{2}}(\gamma(-\delta),\gamma(\delta)))}{\mathfrak{s}_{\kappa}(d_{N,\varepsilon,f}(\gamma(-\delta),\gamma(\delta)))}&=1+\frac{(1-\eps)g(\nabla f,v)}{n-1}\delta+\frac{\kappa}{2}e^{\frac{-4(1-\eps)f(x)}{n-1}}\delta^2+O(\delta^3).
\end{align*}
From $(1+s)^{a}=1+a s+(a(a-1)/2)s^2+O(s^3)$,
we conclude \eqref{eq:app2}.

We next sketch the proof of \eqref{eq:app6}.
Combining \eqref{eq:app2}, \eqref{eq:app3}, \eqref{eq:app4},
we have
\begin{align*}
&\quad \,\, e^{-f(\gamma(-\delta))+f(x)}\,(1+\theta_{\eps}\,\delta)^{n}\,\beta_{\kappa,N,\varepsilon,f,\frac{1}{2}}(\gamma(\delta),\gamma(-\delta)) \\
&= 1+\left\{n\theta_{\eps}-\left(c^{-1}\frac{1-\eps}{n-1}-1\right)g(\nabla f,v)\right\}\delta -n\theta_{\eps} \left(c^{-1}\frac{1-\eps}{n-1}-1\right)g(\nabla f,v)\delta^2\\ 
&+\left[ \left\{ \left( 1-c^{-1}\frac{1-\eps}{n-1}  \right)^2-c^{-1}\left(\frac{1-\eps}{n-1}\right)^2      \right\}g(\nabla f,v)^{2}\right]\frac{\delta^{2}}{2}\\ 
&+\left(c^{-1}\,\kappa\,e^{\frac{-4(1-\eps)f(x)}{n-1}}-\nabla^2 f(v,v)+n(n-1) \theta^2_{\eps}\right)\frac{\delta^{2}}{2}+O(\delta^3),\\
&\quad \,\, e^{-f(\gamma(\delta))+f(x)}\,(1-\theta_{\eps}\,\delta)^{n}\,\beta_{\kappa,N,\varepsilon,f,\frac{1}{2}}(\gamma(-\delta),\gamma(\delta))\\
&= 1-\left\{n\theta_{\eps}-\left(c^{-1}\frac{1-\eps}{n-1}-1\right)g(\nabla f,v)\right\}\delta-n\theta_{\eps} \left(c^{-1}\frac{1-\eps}{n-1}-1\right)g(\nabla f,v)\delta^2\\ 
&+\left[ \left\{ \left( 1-c^{-1}\frac{1-\eps}{n-1}  \right)^2-c^{-1}\left(\frac{1-\eps}{n-1}\right)^2      \right\}g(\nabla f,v)^{2}\right]\frac{\delta^{2}}{2}\\
&+\left(c^{-1}\,\kappa\,e^{\frac{-4(1-\eps)f(x)}{n-1}}-\nabla^2 f(v,v)+n(n-1)\theta^2_{\eps}\right)\frac{\delta^{2}}{2}+O(\delta^3).
\end{align*}
This implies
\begin{align}\notag
&\left(e^{-f(\gamma(-\delta))+f(x)}\,(1+\theta_{\eps}\,\delta)^{n}\beta_{\kappa,N,\varepsilon,f,\frac{1}{2}}(\gamma(\delta),\gamma(-\delta))\right)^{\frac{c}{c+1}}\\ \notag
&= 1+\frac{c}{c+1}  \left\{n\theta_{\eps}- \left(c^{-1}\alpha-1\right)g(\nabla f,v)\right\}\delta \\ \notag
&+\frac{c}{c+1}\left(c^{-1}\,\kappa\,e^{\frac{-4(1-\epsilon)f(x)}{n-1}}-\nabla^2 f(v,v)-\frac{F(\theta_{\eps})}{c+1}\right)\frac{\delta^{2}}{2}+O(\delta^3),\\ \notag
&\left(e^{-f(\gamma(\delta))+f(x)}\,(1-\theta_{\eps}\,\delta)^{n}\beta_{\kappa,N,\varepsilon,f,\frac{1}{2}}(\gamma(-\delta),\gamma(\delta))\right)^{\frac{c}{c+1}}\\ \notag
&= 1-\frac{c}{c+1}  \left\{n\theta_{\eps}- \left(c^{-1}\alpha-1\right)g(\nabla f,v)\right\}\delta \\ \notag
&+\frac{c}{c+1}\left(c^{-1}\,\kappa\,e^{\frac{-4(1-\epsilon)f(x)}{n-1}}-\nabla^2 f(v,v)-\frac{F(\theta_{\eps})}{c+1}\right)\frac{\delta^{2}}{2}+O(\delta^3),
\end{align}
where $F(\theta)$ is defined as \eqref{eq:app5}.
In particular,
\begin{align}\notag
&\frac{1}{2}\left(e^{-f(\gamma(-\delta))+f(x)}\,(1+\theta_{\eps}\,\delta)^{n}\,\beta_{\kappa,N,\varepsilon,f,\frac{1}{2}}(\gamma(\delta),\gamma(-\delta))\right)^{\frac{c}{c+1}}\\ \notag
&+\frac{1}{2}\,\left(e^{-f(\gamma(\delta))+f(x)}\,(1-\theta_{\eps}\,\delta)^{n}\,\beta_{\kappa,N,\varepsilon,f,\frac{1}{2}}(\gamma(-\delta),\gamma(\delta))\right)^{\frac{c}{c+1}}\\ \notag
&= 1+\frac{c}{c+1}\left(c^{-1}\,\kappa\,e^{\frac{-4(1-\epsilon)f(x)}{n-1}}-\nabla^2 f(v,v)-\frac{F(\theta_{\eps})}{c+1}\right)\frac{\delta^{2}}{2}+O(\delta^3).
\end{align}
Substituting this equation into \eqref{eq:BM1}, we arrive at \eqref{eq:app6}.

\bigskip

\noindent
\emph{Acknowledgment.}
The authors are grateful to the anonymous referee for valuable comments.

\providecommand{\bysame}{\leavevmode\hbox to3em{\hrulefill}\thinspace}

\end{document}